\documentclass[letterpaper, 12pt]{iopart}
\expandafter\let\csname equation*\endcsname\relax
\expandafter\let\csname endequation*\endcsname\relax
\usepackage{amsmath, amssymb}

\usepackage{etex}
\usepackage{mathrsfs}
\usepackage{microtype}
\usepackage{booktabs}
\usepackage{tabularx}
\usepackage{pifont}
\usepackage{overpic}
\usepackage[square, comma, sort&compress, numbers]{natbib}
\usepackage[center,sc]{caption}

\usepackage{subfigure}

\usepackage{fancyvrb}

\usepackage[pdfauthor={Philippe H. Trinh, S. Jonathan Chapman},pdftitle={Exponential asymptotics and coalescing corners}]{hyperref}

\def\Xint#1{\mathchoice
   {\XXint\displaystyle\textstyle{#1}}%
   {\XXint\textstyle\scriptstyle{#1}}%
   {\XXint\scriptstyle\scriptscriptstyle{#1}}%
   {\XXint\scriptscriptstyle\scriptscriptstyle{#1}}%
   \!\int}
\def\XXint#1#2#3{{\setbox0=\hbox{$#1{#2#3}{\int}$}
     \vcenter{\hbox{$#2#3$}}\kern-.5\wd0}}

\def\dashint{\Xint-}

\renewcommand*{\Re}{\operatorname{Re}} 
\renewcommand*{\Im}{\operatorname{Im}} 
\newcommand*{\Ei}{\operatorname{Ei}}
\renewcommand*{\L}{\operatorname{\mathfrak{L}}}
\newcommand*{\Nop}{\operatorname{\mathfrak{N}}}

\newcommand*{\ep}{\epsilon}
\newcommand*{\N}{\mathcal{N}}
\newcommand*{\im}{\mathrm{i}}
\newcommand*{\St}{\mathcal{S}}
\newcommand*{\Omcc}{\Omega^\text{cc}}
\newcommand*{\phicc}{\phi_\text{exp}^\text{cc}}
\newcommand*{\phisingle}{\phi_\text{exp}^\text{single}}

\newcommand{\de}{\textrm{d}}
\newcommand{\dd}[2]{\frac{\de#1}{\de#2}}

\newcommand{\Oh}{\mathcal{O}}

\renewcommand{\e}{\mathrm{e}}
\newcommand{\eg}{\emph{e.g.} }
\renewcommand*{\etal}{\emph{et al.} }
\newcommand*{\sqep}{\ep^{1/2}}

\begin{document}

\title[Exponential asymptotics with coalescing singularities]{Exponential asymptotics with coalescing singularities}

\author{
Philippe H. Trinh and S. Jonathan Chapman}

\address{Oxford Centre for Industrial and Applied Mathematics, Mathematical Institute, University of Oxford}
\ead{trinh@maths.ox.ac.uk, chapman@maths.ox.ac.uk}
\begin{abstract}
Problems in exponential asymptotics are typically characterized by divergence of the associated asymptotic expansion in the form of a factorial divided by a power. In this paper, we demonstrate that in certain classes of problems that involve coalescing singularities, a more general type of exponential-over-power divergence must be used. As a model example, we study the water waves produced by flow past an obstruction such as a surface-piercing ship. In the low speed or low Froude limit, the resultant water waves are exponentially small, and their formation is attributed to the singularities in the geometry of the obstruction. However, in cases where the singularities are closely spaced, the usual asymptotic theory fails. We present both a general asymptotic framework for handling such problems of coalescing singularities, and provide numerical and asymptotic results for particular examples. 
\end{abstract}


\section{Introduction}

Many problems in exponential asymptotics involve the analysis of singularly perturbed differential equations where the associated solution is expressed as a divergent asymptotic expansion. It has been noted by authors such as \citeauthor{dingle_1973} \cite{dingle_1973} and \citeauthor{berry_1991} \cite{berry_1991} that in many cases, the divergence of the sequence occurs in the form of a factorial divided by a power. In this paper, we use a model problem from the theory of water waves and ship hydrodynamics to demonstrate how in certain classes of problems, a more general form of divergence must be used in order to perform the exponential asymptotic analysis. This behaviour is expected to occur in a wide range of singular perturbation problems characterized by multiple singularities coalescing in the limit the small parameter tends to zero. 

We begin with a brief explanation of exponential asymptotics and factorial-over-power divergence. Consider a re-scaled differential equation for the exponential integral function,
\begin{equation} \label{eiode}
	\epsilon \dd{y}{z} + y = \frac{\epsilon}{z}, \qquad \text{$y \to 0$ as $z \to -\infty$},
\end{equation}
in the limit $\epsilon \to 0$, where $y: \mathbb{C} \mapsto \mathbb{C}$, and for which the exact solution is given by $y(z) = e^{-z/\epsilon}\Ei(z/\epsilon)$. Begin by assuming that $\Im(z) < 0$ in \eqref{eiode}. The usual approach is to expand the solution as a series 
\begin{equation} \label{ybase}
	y(z) = \sum_{n=1}^\infty \epsilon^n y_n,
\end{equation}
through which it is found that $y_n = (n-1)!/z^n$. However, writing \eqref{eiode} as an integral for $y(z)$, we find that if $z$ is analytically continued along a path that crosses the positive real axis, an exponentially small term switches-on across $\Re(z) \geq 0$. Once this occurs, the asymptotic form of \eqref{ybase} is modified to 
\begin{equation}
	y(z) = \sum_{n=1}^\infty \epsilon^n y_n + 2\pi \im \e^{-z/\epsilon}.
\end{equation}
The unusual process by which exponentially small terms can suddenly appear or disappear in an asymptotic expansion is known as the \emph{Stokes Phenomenon} (\citeauthor{berry_1989} \cite{berry_1989}, \citeauthor{meyer_1989} \cite{meyer_1989}, \citeauthor{olde_1995} \cite{olde_1995}). In the case of \eqref{eiode}, the appearance of the exponentially small terms can be understood through a variety of techniques ranging from the method of steepest descents applied to the integral expression (\citeauthor{bleistein_1975} \cite{bleistein_1975}), to Borel transforms (\citeauthor{dingle_1973} \cite{dingle_1973}), optimal truncation of the regular expansion (\citeauthor{chapman_1998} \cite{chapman_1998}), or techniques of series acceleration (\citeauthor{baker_1975} \cite{baker_1975}). 

The key feature in such problems exhibiting the Stokes Phenomenon is the divergence of the na\"{i}ve asymptotic expansion \eqref{ybase}. A generic quality of singular perturbation problems is that once \eqref{ybase} has been truncated at $n = \mathcal{N}$, say, an equation for the leading-order remainder, $R_\mathcal{N}$, will be of the form
\begin{equation}
	\L(R_\mathcal{N}; \epsilon) 
    \sim \epsilon^{\mathcal{N}} y_\mathcal{N},
\end{equation}
where $\L$ is a linear differential operator. However, in the limit $\epsilon \to 0$, the optimal truncation point, $\mathcal{N} \to \infty$, and hence the exponentially small remainder is related to the divergent behaviour of $y_\N$. It has been noted, principally by \citeauthor{dingle_1973} \cite{dingle_1973}, that in singularly perturbed problems, the divergence of the asymptotic series takes the form
\begin{equation}
	y_n \sim \frac{\Gamma(n + a)}{(\text{variable})^{n}},
\end{equation}
as $n\to \infty$, and this `factorial-over-power' divergence of the late terms is understood as the consequence of expanding a function with isolated singularities in the complex plane, and is a result also known as Darboux's Theorem (Henrici \cite[p.447]{henrici_1977}).

In this paper, we demonstrate that for a certain problems that involves coalescing singularities, the divergence of the late terms may take the alternative form of an exponential-over-power, namely
\begin{equation} \label{exppower}
	y_n \sim \frac{\Gamma(n + a)}{(\text{variable})^{n}}
	\exp\biggl[ \sum_j \text{(variable)} n^{a_j}\biggr],
\end{equation}
where $0 < a_j < 1$.

In particular, the theory we describe is applicable for the study of certain singular differential equations of the form 
\begin{equation} \label{Nop}
	\Nop(z, y; \epsilon, a_1(\epsilon), a_2(\epsilon), \ldots) = 0,
\end{equation}
where $\Nop$ is a nonlinear differential operator on $y: \mathbb{C} \mapsto \mathbb{C}$, $z\in\mathbb{C}$, and it is assumed that the asymptotic expansion of $y$ in the limit $\ep \to 0$ contains singularities at points $z = -a_k$, in the complex plane. The asymptotic expansion will diverge based on the presence of these singularities, and Stokes lines will necessitate the switching-on of exponentially small terms. However, there exists a distinguished limit whereby the singularities may coalesce (\emph{e.g.} $a_i \to a_j$ for distinct $i$, $j$), and the study of this distinguished limit is the subject of this paper. 

Although our aim is to present a general methodology for such nonlinear differential equations, the exponential-over-power form of \eqref{exppower} commonly arises in linear differential and difference equations, so we shall begin in Sec.~\ref{sec:linear} by illustrating these connections through a simple example.

We continue in Sec.~\ref{mathform} with a presentation of a model problem motivated by previous work on the application of exponential asymptotics to the study of surface waves produced by flow over a submerged object (Chapman and Vanden-Broeck \cite{chapman_2002, chapman_2006}, Trinh and Chapman \cite{trinh_2011, trinh_2014, trinh_2013, trinh_2013a}). In the limit of low Froude numbers (representing the balance between inertial and gravitational forces), potential flow past an obstruction, such as a step in a channel or a surface-piercing ship, will produce free-surface waves that are exponentially small in the Froude number. These low-Froude water wave problems present a useful setting for studies on exponential asymptotics, as the waves can be directly observed in numerical and experimental settings, and concepts such as Stokes lines and the Stokes Phenomenon share a correspondence with the physical fluid domain.

In Sec.~{\ref{sec:wellsep}}, we review the application of exponential asymptotics to study the case where the singularities are well separated. The techniques we apply are based on the use of a factorial over power ansatz to capture the divergence of the asymptotic expansions, then optimal truncation and Stokes line smoothing to relate the late-order terms to the exponentially small waves (see for example, papers by \citeauthor{olde_1995} \cite{olde_1995}, \citeauthor{chapman_1998} \cite{chapman_1998}, and \citeauthor{trinh_2010} \cite{trinh_2010}). For more details on standard techniques in exponential asymptotics, including other approaches, we refer the reader to the tutorials and reviews by Boyd (\citeyear[Chap. 4]{boyd_1998}), \citeauthor{olde_1999} \cite{olde_1999}, \citeauthor{costin_2008} \cite{costin_2008}, and \citeauthor{grimshaw_2010} \cite{grimshaw_2010}.

The problem with coalescing singularities is then studied in Secs.~{\ref{sec:cc13}} to \ref{sec:nhull_cc1616}, for a particular case, while the most general methodology is presented in \ref{sec:general}.

\section{A toy linear differential equation} \label{sec:linear}

The exponential-over-power expression of the sort in \eqref{exppower} is a familiar sight to those who are well aquainted with the study of linear difference equations or linear differential equations near an irregular singular point. In particular, there is a connection between our exponential asymptotic analysis, the exponential-over-power expression \eqref{exppower}, and the asymptotic theory of linear difference equations that was first established in the classic papers of \citeauthor{birkhoff_1930} \cite{birkhoff_1930} and \citeauthor{birkhoff_1933} \cite{birkhoff_1933} (refer to \citeauthor{wimp_1985} \cite{wimp_1985} for a more readeable treatment). 

Consider as a toy example the differential equation
\begin{equation} \label{toy1}
\ep \left[1 + \frac{\ep^\alpha}{x^{1/2}}\right] y'(x) + y(x) = \frac{1}{x^{1/2}} + \frac{\ep^{1/2}}{x}.    
\end{equation}
where we first set $\alpha = 1$. The solution can be written as a regular expansion of the form $y = \sum \ep^{n/2} y_n$, and we see that the first two orders, $y_0 = 1/x^{1/2}$ and $y_1 = 1/x$ are singular at $x = 0$. Since all subsequent orders depend on derivatives of these two terms, it is expected that the series is divergent. At $\Oh(\ep^{n/2})$, we have
\begin{equation}
y_n = -y_{n-2}' - \frac{y_{n-4}'}{x^{1/2}},
\end{equation}
and it can be verified that in the limit $n \to \infty$, the divergence is captured by the standard factorial-over-power ansatz,
\begin{equation} \label{forman}
    y_n \sim \frac{A(x) \Gamma(\frac{n}{2})}{x^{n/2}}.
\end{equation}
where $A(x) = A_0 \e^{2\sqrt{x}}$ and $A_0$ is constant. 

Now consider a modification to the factor multiplying $y'$ in \eqref{toy1}. If we set $\alpha = 1/2$, then the $\Oh(\ep^{n/2})$ equation changes to 
\begin{equation}
y_n = -y_{n-2}' - \frac{y_{n-3}'}{x^{1/2}},
\end{equation}
and we find that the ansatz \eqref{forman} no longer works because the $y_{n-3}$ term contributes before the prefactor $A(x)$ can be determined. The required ansatz is instead of exponential-over-power form,
\begin{equation} \label{yntoy}
y_n \sim \frac{B(x)\Gamma(\frac{n}{2}) \e^{\sqrt{2n}}}{x^{n/2}},
 \end{equation} 
where $B(x) = B_0/x$ and $B_0$ is constant.

We may also relate the differential equation \eqref{toy1} to a linear homogeneous difference equation. If we re-scale near the singularity, with $x = \ep z$ and $y(x) = Y(z)/\epsilon^{1/2}$, then the differential equation becomes 
\begin{equation} \label{toy2}
\left( 1 + \frac{1}{z^{1/2}}\right) Y'(z) + Y(z) = \frac{1}{z^{1/2}} + \frac{1}{z}.
\end{equation}
Then, seeking an expansion of the form $Y = \sum A_n/z^{n/2}$, valid in the limit $z \to \infty$, gives the recurrence relation, $A_1 = 1$, $A_2 = 1$, and in general for $n \geq 3$,
\begin{gather} \label{toyrec}
A_n = \left(\frac{n}{2} - 1\right)A_{n-2} + \left(\frac{n}{2}-\frac{3}{2}\right) A_{n-3}
\end{gather}
In the limit $n \to \infty$, the divergence of $A_n$ follows $A_n \sim \Lambda\Gamma(n/2)\e^{\sqrt{2n}}$ for constant $\Lambda$. This form indeed matches the outer-to-inner limit of \eqref{yntoy} as the singularity is approached. 
	
In fact, we shall find in Sec.~\ref{sec:cc13} that the toy model \eqref{toy2} is related to the leading-order inner solution of the main nonlinear problem \eqref{phi_ode} of this paper [compare the recurrence relation \eqref{toyrec} with \eqref{recAn}]. For the main problem of interest, the coalescence of singularities will cause the respective inner equation to resemble \eqref{toy2}, provided that the rate of coalescence is chosen appropriately in the limit $\ep \to 0$.

\section{A model for the ship-wave problem} \label{mathform}

Although the main contents of this paper can be appreciated without understanding the physical context of the differential equations, it is still helpful to understand from where the model arose, and the relationships between mathematical and physical theories.  

The governing equations for two-dimensional, steady, incompressible, irrotational, and inviscid flow past a submerged or surface-piercing object in the presence of gravity involves: (i) the solution of Laplace's equation  for the fluid potential, $\phi$; (ii) kinematic conditions, on all fluid and solid surfaces; and (iii) a dynamic boundary condition (Bernoulli's equation) for the free surface. In the low speed or low Froude limit, a small non-dimensional parameter, $\epsilon$, representing the balance between inertial and gravitational forces, can be introduced.

We explain in \ref{app:realfluid} that the search for exponentially small free-surface waves in this system can be modeled using the complex initial value problem,
\begin{gather} \label{phi_ode} 
\Bigl[ \phi - q_s^2 \Bigr] - \im\epsilon q_s \phi \dd{\phi}{w} = 0, 
\end{gather}
with $\phi(0) = 0$, where $w \in \mathbb{R}$, and $\phi: \mathbb{C} \to \mathbb{C}$. Physically, $w$ corresponds to the complex potential and $\phi$ is related to the square of the fluid velocity. In this paper, we are primarily interested in studying the case where the forcing function, $q_s$, is given by 
\begin{equation} \label{qs_2hull}
q_s(w) = \frac{w^{\sigma_1+\sigma_2}}{(w+a_1)^{\sigma_1} (w+a_2)^{\sigma_2}}, 
\end{equation}
where $0 < a_2 < a_1$, $a_1 + a_2 = 1$, and $0 < \sigma_1, \sigma_2 < 1$. We note that in the limit $\epsilon \to 0$, $\phi \sim q_s^2$. In terms of the fluid mechanics, this $q_s$ function describes the leading-order speed for low-Froude flow past a two-cornered ship with corners located at $w = -a_1$ and $w = -a_2$ in the potential plane, and with divergent corner-angles $\pi \sigma_1$ and $\pi \sigma_2$. 

\begin{figure}[htb]\centering
\includegraphics{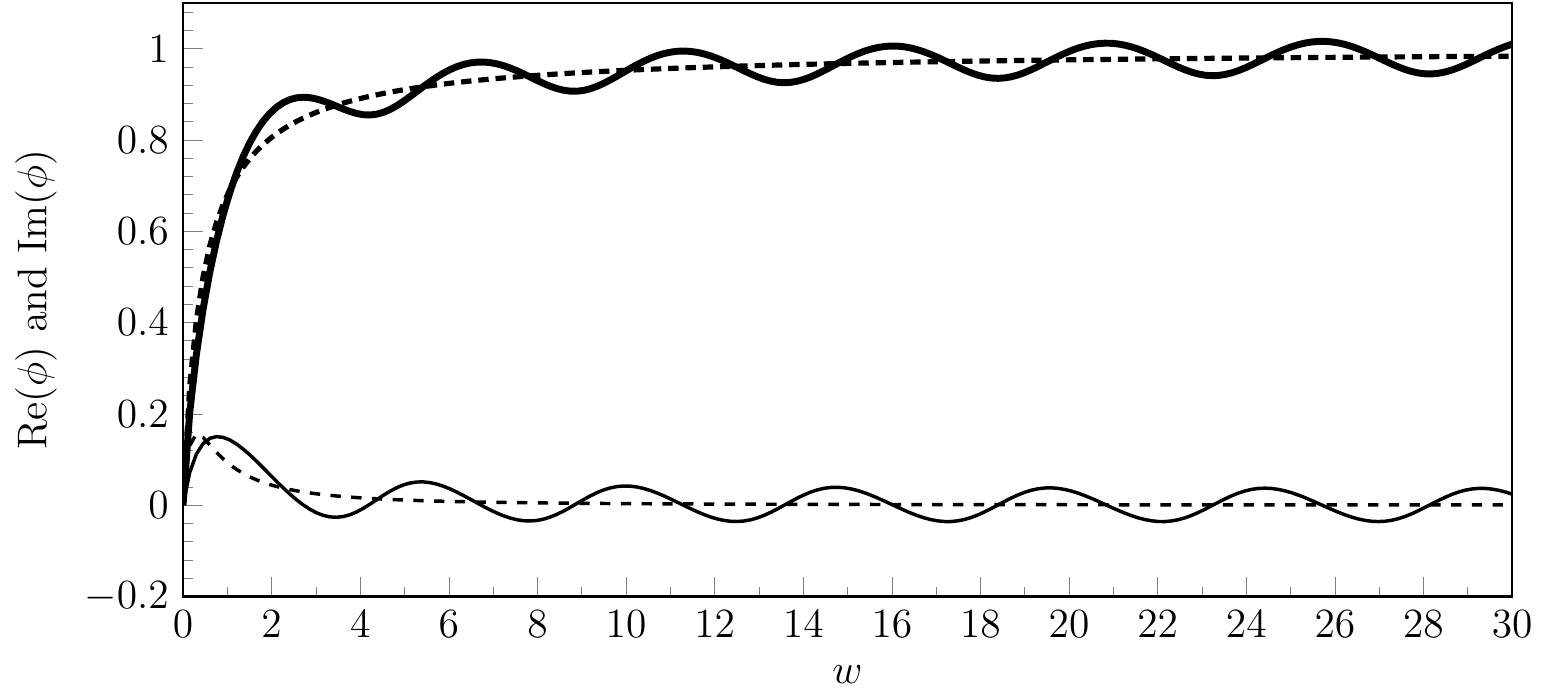}
\caption{Numerical solutions $\Re(\phi)$ (thick) and $\Im(\phi)$ (thin) for $\epsilon = 0.8$, $a_1 = 0.8$, $a_2 = 0.2$, $\sigma_1 = \sigma_2 = 1/4$. The dashed lines correspond to $q_s^2$ (thick) and $2\im q_s^4 q_s'$ (thin). \label{phiprofile}} 
\end{figure}

A numerical solution of the differential equation for $\epsilon = 0.8$ is shown in Figure \ref{phiprofile}. Here, $a_1 = 0.8$, $a_2 = 0.2$, and $\sigma_1 = \sigma_2 = 1/4$. We note that because of the singularity as $w \to 0$, the numerical solution is solved beginning near $w = 0$ ($w_0 = 10^{-5}$ in the figure) subject to the initial condition of $\phi= q_s^2 + 2\epsilon \im q_s^4 q_s'$ evaluated at $w = w_0$ (this expression is derived from the asymptotic expansion as $\epsilon \to 0$ and is covered later).

Our goal in this paper is to derive the form of the waves appearing in Figure \ref{phiprofile} for the particular case when the two singularities, $a_1$ and $a_2$, in \eqref{qs_2hull} tend to one another as $\epsilon \to 0$. In the limit that $a_1 \to a_2 = a$, we intuitively expect that the forcing function can be replaced by the single-singularity forcing,
\begin{equation} \label{qs_one}
	q_s(w) = \left(\frac{w}{w+a}\right)^{\sigma},
\end{equation}
with $a > 0$ and $\sigma = \sigma_1 + \sigma_2$. However, as we shall see, the replacement of \eqref{qs_2hull} by \eqref{qs_one} is non-trivial, and several challenging aspects emerge when considering the asymptotic analysis in the distinguished limit of merging singularities. 

Although the full problem in \eqref{governingreal} can be studied using our methods, there are two principle reasons why we prefer to study the differential equation \eqref{phi_ode}. Firstly, the simplification eases the algebraic complexity of the asymptotic analysis, while still preserving all the features of the theory we wish to present (in relation to the coalescence of singularities). Secondly, accurate numerical verification of the asymptotic analysis demands numerical precision of five or six digits of accuracy---otherwise, the fine effects of adjusting the ship's geometry are easily missed; this precision can only be easily achieved for the simpler problem, which does not require the computation of the Cauchy Principal integral of the full governing equations. 

In particular, we note the qualitative similarities between solutions of our reduced problem in Figure~\ref{phiprofile}, and the full boundary-integral solutions displayed in Figure~7 of \cite{trinh_2011} and Figure~4 of \cite{trinh_2014}. A similar idea of reducing the boundary integral was used by Tuck \cite{tuck_1991}, who realized that the integral does not play a significant role in seeking free-surface waves in low-Froude problems.

\section{A review of the case of well-separated singularities} \label{sec:wellsep}

The theory for the nonlinear equations of low-speed ship flows with well-separated singularities is presented in \cite{trinh_2014}. Here, although we study a slightly different problem \eqref{phi_ode}, we shall review the main procedure for applying exponential asymptotics, with particular emphasis on the breakdown when the singularities are closely spaced. 

\subsection{Step 1: Characterize divergence of late terms} \label{sec:divergence}

\noindent We begin as usual by applying a regular asymptotic expansion to \eqref{phi_ode}
\begin{equation} \label{eq:phiseries}
 \phi = \sum_{n=0}^\infty \epsilon^n \phi_n, 
\end{equation}

\noindent valid in the limit $\epsilon \to 0$. The leading-order solution is related to the rigid-wall flow of \eqref{qs_one}, and given by 
\begin{equation} \label{phi0_2hull}
	\phi_0 = q_s^2 = \left(\frac{w^{\sigma_1+\sigma_2}}{(w+a_1)^{\sigma_1} (w+a_2)^{\sigma_2}}\right)^2.
\end{equation}
We note the presence of singularities at $w = -a_1$ and $w = -a_2$, which are related to the physical corners of the ship. At $\Oh(\epsilon^n)$, we obtain
\begin{equation} \label{phin}
	\phi_n = \im \epsilon q_s \Bigl[ \Bigl(\phi_0 \phi_{n-1}' + \phi_1 \phi_{n-2}' + \ldots \Bigr) + \Bigl( \phi_{n-1} \phi_0' + \ldots\Bigr)\Bigr].
\end{equation}
The calculation of $\phi_n$ at each order relies upon the differentiation of $\phi_{n-1}$, and thus at each order of the asymptotic procedure, we increase the power of the singular term from the previous order. In the limit $n \to \infty$, the late terms diverge as a factorial over power of the form,
\begin{equation} \label{ansatz}
	\phi_n \sim \frac{P_1(w) \Gamma(n + \gamma_1)}{[\chi_1(w)]^{n+\gamma_1}} + \frac{P_2(w) \Gamma(n + \gamma_2)}{[\chi_2(w)]^{n+\gamma_2}}.
\end{equation}
where $\gamma_k$ are constant, and $\chi = \chi_k$ are known as the \emph{singulants}, with $\chi_1(-a_1) = 0 = \chi_2(-a_2)$ from the two singularities in the leading-order equation \eqref{phi0_2hull}. Once the ansatz \eqref{ansatz} is substituted into \eqref{phin}, we obtain from the leading-order contribution as $n \to \infty$ the expressions for $\chi_1$ and $\chi_2$, given by 
\begin{equation} \label{chik}
\chi_k(w) = \int_{-a_k}^w \frac{\im}{q_s^3(\varphi)} \, \de{\varphi}. 
\end{equation}

\begin{figure} \centering
\includegraphics{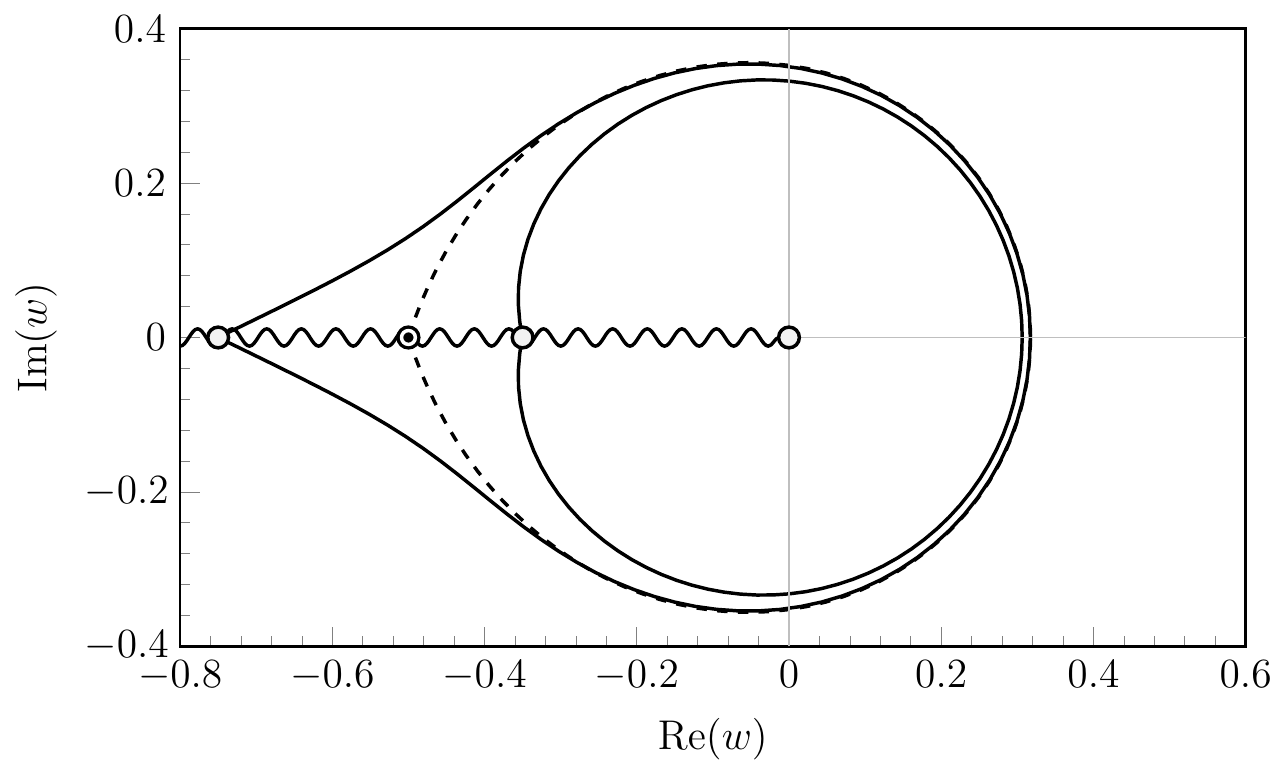}
\caption{(Solid) Stokes lines, $\Im(\chi_k) = 0$ and $\Re(\chi_k) \geq 0$, using \eqref{qs_2hull} with $a_1 = 0.75$, $a_2 = 0.35$ and $\sigma_1 = \sigma_2 = 1/4$. (Dashed) Stokes line for \eqref{qs_one} for $a = 0.5$ and $\sigma = 1/2$. The branch cuts for each of the critical points in $q_s$ is taken along the negative real axis and shown as a snaking curve. \label{fig:stokes}}
\end{figure}

Using the expression of $\chi_k$, Stokes lines can be traced from each of the two singularities, across which the Stokes Phenomenon necessitates the switching-on of waves.  From Dingle \cite{dingle_1973}, these special lines are given at the points, $w\in\mathbb{C}$, where
\begin{equation}
\Im[\chi_k(w)] = 0 \text{\quad and \quad} \Re[\chi_k(w)] \geq 0.	
\end{equation}
The Stokes lines are computed by numerically integrating \eqref{chik} and are shown in Figure \ref{fig:stokes} for a particular choice of $a_k$ and $\sigma_k$. In the figure, we observe Stokes lines emerging from $w = -a_1$ and $w = -a_2$, and intersecting the positive real $w$-axis. It is expected that across these two points of intersection, an exponential switches on. Also shown (dashed) in the figure is the Stokes line that corresponds to the one-singularity function \eqref{qs_one}.

At the subsequent order in \eqref{phin} as $n\to\infty$, using the ansatz \eqref{ansatz}, we obtain
\begin{equation} \label{Pk}
P_k(w) = \frac{\Lambda_k}{q_s^4(w)},
\end{equation}
where $\Lambda_k$ is a constant prefactor that must be calculated from matching the outer asymptotic expansion \eqref{eq:phiseries} with the nonlinear solution near the singularities, $w = -a_k$. Finally, the value of $\gamma_k$ can be derived by requiring the singular behaviour of $\phi_n$ in \eqref{ansatz} to match $\phi_0$ in \eqref{phi0_2hull} when $n = 0$. This gives a value of 
\begin{equation} \label{gammak}
	\gamma_k = \frac{6\sigma_k}{1 + 3 \sigma_k}.
\end{equation}

\subsection{Step 2: Match with the solution near the singularity} \label{sec:match}

In order to determine the unknown pre-factors $\Lambda_k$ that characterize the divergence of the late terms in \eqref{ansatz}, we must rescale $w$ and $\phi$ near the singularities, $w = -a_k$, and match to the outer solution found in the previous section. The size of this inner region can be derived by observing where the breakdown in the outer expansion (\ref{eq:phiseries}) first occurs, \emph{i.e.} where $\epsilon q_1 \sim q_0$. This inner region is then delimited by
\begin{equation} \label{innerscaling}
 \phi = \Oh( q_s^2) = \Oh\Bigl( \bigl(w + a_k \bigr)^{-2\sigma_k}\Bigr) \quad \text{and} \quad 
 w+a_k \sim \mathcal{O}\left(\epsilon^{\frac{1}{1+3\sigma_k}}\right).
\end{equation}
Once $\phi$ and $w$ have been rescaled with \eqref{innerscaling} in mind, then a recurrence relation can be developed for the inner solution, which is then numerically solved (if required). This inner solution is then matched with the outer solution \eqref{eq:phiseries} in order to determine $\Lambda_k$. To be specific, we note that as $w \to -a_k$, 
\begin{equation}
	q_s \sim c_k (w + a_k)^{-\sigma_k},
\end{equation}
where the complex constant $c_k$ is computed from \eqref{qs_2hull}. Within the inner region, the solution, $\phi$, can be expanded as an infinite series,
\begin{equation} \label{phiinner}
\phi = c_k^2 (w+a_k)^{-2\sigma_k} \sum_{n=0}^\infty \frac{A_n}{z^n},
\end{equation}
where $z = X(w+a)^{1 + 3\sigma_k}/\epsilon$, with $X = \im/[c_k^3(1 + 3\sigma_k)]$. We then find the coefficients, $A_n$, by solving the recurrence relation 
\begin{equation}
A_0 = 1, \qquad A_n = \sum_{m=0}^{n-1} \left(m + \frac{2\sigma_k}{1 + 3\sigma_k}\right) A_m A_{n-m-1}, \quad \text{for $n > 1$}. 	
\end{equation}
The divergence of the coefficients, $A_n$, is described through the constant,
\begin{equation}
\Omega(\sigma_k) = \lim_{n\to\infty} \frac{A_n}{\Gamma(n+\gamma_k)},
\end{equation}
which can be computed numerically. Finally, matching the outer series \eqref{eq:phiseries} with the inner series \eqref{phiinner} gives a value for the unknown prefactor that appears in the late terms:  
\begin{equation}
\Lambda_k = \frac{c_k^{6 - 3\gamma_k} \e^{\frac{\pi \im}{2} \gamma_k}}{(1 + 3\sigma_k)^{\gamma_k}} \Omega(\sigma_k). \end{equation}
This completes the determination of all components of the late terms in \eqref{ansatz}. 

\subsection{Step 3: Optimally truncate and smooth the Stokes line}

To derive the form of the exponentials that appear whenever a Stokes Line
intersects the free-surface, we optimally truncate the asymptotic
expansions (\ref{eq:phiseries}), and examine the remainder as the Stokes line is crossed. We let
\begin{equation} \label{eq:nhull_opttrun}
 \phi = \sum_{n=0}^{\mathcal{N}-1} \epsilon^n \phi_n + R_\mathcal{N}.
\end{equation}

\noindent When $\mathcal{N}$ is chosen to be the optimal truncation point, the remainder $R_\mathcal{N}$ is found to be exponentially small; by re-scaling near the Stokes line, it can be shown that a wave of the following form switches on:
\begin{equation} \label{2hullphiexp} 
\phi_\text{exp, k} \sim -\frac{2\pi \im}{\epsilon^{\gamma_k}} P_k \exp\left[{-\frac{\chi_k}{\epsilon}}\right].
\end{equation}
where $k = 1$, $2$, and note that such a contribution is only included if the associated Stokes line from $w= -a_k$ intersects the positive real $w$-axis (the negative sign is due to a switching moving in the direction of positive $w$ across the Stokes line). 
%

Since $P_k = \Lambda_k/q_0^4(w) \to \Lambda_k$ as $w \to \infty$ by \eqref{Pk}, then the form of the far-field waves are
\begin{equation} \label{phiexp_two}
 \Re(\phi_{\text{exp}, k}) \sim -\frac{2\pi \Omega(\sigma_k)}{\epsilon^{\gamma_k}} \left[\frac{|c_k|^{6-3\gamma_k}}{(1+3\sigma_k)^{\gamma_k}} \right]
\exp \Biggl[-\frac{\Re(\chi_k)}{\epsilon}\Biggr] \sin \Biggl(
-\frac{w}{\epsilon} + \Psi_k \Biggr) 
\end{equation}

\noindent where the phase shift is
\begin{equation}
\Psi_k = 
-\frac{1}{\epsilon} \Im\left[\int_{-a_k}^{-a_1} \dd{\chi}{\varphi} \ d\varphi \right]
+ (6-3\gamma_k) \text{Arg}(c_k) + \frac{\pi}{2}\gamma_k.
\end{equation}

\begin{figure}[htb] \centering
\includegraphics{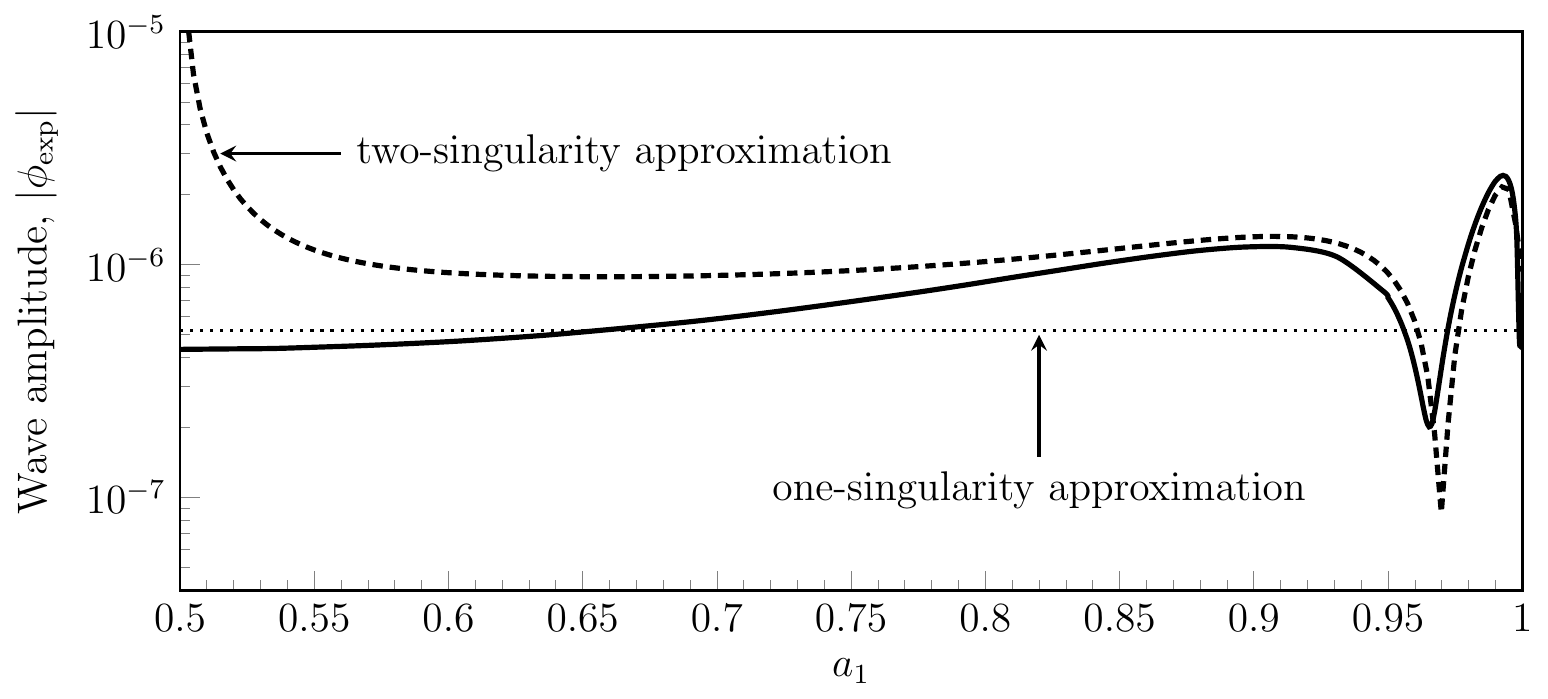}
\caption{The numerical solution (solid) for the $\sigma_1 = \sigma_2 = 1/4$ forcing is plotted against the well-separated asymptotic approximation at $\epsilon = 0.15$ and $a_1 + a_2 = 1$ (dashed). The two-singularity approximation very accurately predicts the solution when the two singularities are well spaced, but is singular when $a_1, a_2 \approx 0.5$ near the left. The dotted line indicates the one-singularity approximation for $\sigma = 1/3$ amd $a = 0.5$. \label{fig:cc2525}}
\end{figure}

We also note that if we substitute the single-singularity function for $q_s$ \eqref{qs_one} into the differential equation \eqref{phi_ode}, then the same exponential form \eqref{phiexp_two} applies. The key observation is that within the outer region, where the ansatz \eqref{ansatz} plays an important role, the precise nature of the $q_s$ function is unimportant. The determination of the prefactor ansatz, $P_k$, requires matching with an inner solution, where $q_s \sim c_k (w + a_k)^{-\sigma_k}$. Since the single-singularity $q_s$ contains the same local behaviour as the double-singularity form, the exponential prediction \eqref{phiexp_two} follows identically with $c = c_k$, $\sigma = \sigma_k$, and so on. 
∂
The existence of a distinguished limit whereby the two singularities of \eqref{qs_2hull} merge in the limit $\epsilon \to 0$ can be observed by computing numerical solutions of \eqref{phi_ode}. For fixed values of $\epsilon$ and $\sigma_1, \sigma_2$, we set $a_2 = 1 - a_1$ and examine $a_1$ in the range $0.5 \leq a_1 \leq 1$. The wave amplitude for such a computation is shown in Figure \ref{fig:cc2525}. We see that the asymptotic prediction \eqref{phiexp_two} is accurate when the singularities are well separated, but diverges once the singularities approach one another. A uniform approximation would need to smoothly match with the one-singularity approximation at one end, and the (separated) two-singularity analysis at the other. 

\section{Two closely spaced singularities with $\sigma = 1/3$} \label{sec:cc13}

From the well-separated analysis of Sec.~\ref{sec:wellsep}, we observed that the inner region is of the size $\Oh(\epsilon^{1/(1 + 3\sigma_k)})$, and thus as $a_2 \to a_1$, the previous methodology breaks down. We shall then study a problem where the singularities are located a distance $\epsilon^{\ell/m}\beta$ on either side of $w = -a$ (with $a$, $\beta > 0$). The positive integers $\ell$ and $m$ are defined such that
\begin{equation} \label{eq:lmsig}
 \frac{\ell}{m} \equiv \frac{1}{1+3\sigma}
\end{equation}
where $\sigma = \sigma_1 + \sigma_2$ and the ratio $\ell/m$ is irreducible. Now from \eqref{qs_2hull}, we obtain the the forcing function given by
\begin{subequations} \label{qsjoin}
\begin{equation}\label{eq:qspre}
q_s = \frac{w^{\sigma}}{(w+a+\epsilon^{\frac{\ell}{m}}\beta)^{\sigma_1}
(w+a-\epsilon^{\frac{\ell}{m}}\beta)^{\sigma_2}},
\end{equation}
or, as an expansion in $\epsilon$, 
\begin{equation} \label{eq:nhull_cc_qs}
q_s = \left(\frac{w}{w+a}\right)^{\sigma} \sum_{n=0}^\infty
\epsilon^{\frac{\ell n}{m}} \left( \frac{\beta}{w+a}\right)^n f_n, 
\end{equation}
where $f_0 = 1$, $f_1 = \sigma_2 - \sigma_1$, $f_2 = \tfrac{1}{2}[(\sigma_1-\sigma_2)^2 + \sigma_1 + \sigma_2]$, and in general,
\begin{equation} \label{eq:nhull_cc_fn}
f_n = \frac{1}{\Gamma(\sigma_1) \Gamma(\sigma_2)}
\sum_{m=0}^n (-1)^m \frac{\Gamma(\sigma_1+m) \Gamma(\sigma_2+n-m)}
{\Gamma(m+1) \Gamma(n-m+1)}.
\end{equation}
\end{subequations}

Note that the governing equation (\ref{phi_ode}) naturally leads to an expansion in powers of  $\epsilon$, but the fact that $q_s$ is a series in $\epsilon^{\ell/m}$ forces us to expand $\phi$ in more finely spaced powers of $\epsilon^{1/m}$. Our approach, then, is only applicable for cases where $\sigma_1$ and $\sigma_2$ are rational; this is a key requirement in order for the expansions of $\phi$ and $q_s$ to `interleave' properly. Irrational values of $\sigma_1$ and $\sigma_2$ will require more general asymptotic expansions and is beyond the scope of this work (see Appendix A in \cite{trinh_2011} for an example). 

Although it is a straightforward to extend the methodology to handle general rational values of $\sigma_1$ and $\sigma_2$, we will begin by studying the particular case of $\sigma_1 + \sigma_2 = 1/3$ for the forcing $q_s$ in \eqref{eq:nhull_cc_qs}. The general methodology is presented in \ref{sec:general}. 

We set $\ell = 1$ and $m = 2$ in (\ref{eq:lmsig}), and from \eqref{qsjoin}, 
\begin{subequations}
\begin{equation} \label{qs13}
 q_s(w) = \frac{w^{1/3}}
{(w+a+\sqep\beta)^{\sigma_1}(w+a-\sqep\beta)^{\sigma_2}}
= q_0 \sum_{n=0}^\infty \epsilon^{\frac{n}{2}} e_n, 
\end{equation}
with $q_0 = [w/(w+a)]^{1/3}$ and where we have defined the more convenient series coefficients
\begin{equation} \label{en}
e_n = \left(\frac{\beta}{w+a}\right)^n f_n,
\end{equation}
\end{subequations}
with $f_n$ as in (\ref{eq:nhull_cc_fn}).

\subsection{Step 1: Characterize the divergence of the late terms} \label{sec:outercc}

We substitute the regular asymptotic expansion
\begin{equation} \label{phiseries13}
\phi = \sum_{n=0}^\infty \epsilon^{n/2} \phi_n,	
\end{equation}
into the differential equation \eqref{phi_ode}, giving for the low-order terms,
\begin{subequations} \label{earlyorders}
\begin{alignat}{3}
 \Oh(1) 
 \!&:& \qquad 
 \phi_0 &= q_0^2 = \left(\frac{w}{w+a}\right)^{2/3}, \label{phi0} \\
 \Oh(\ep^{1/2}) 
 \!&:& \qquad \phi_1 &= q_0^2 \sum_{k=0}^1 e_k e_{1-k}, \\
 \Oh(\ep) \!&:& \qquad \phi_2 &= q_0^2 \sum_{k=0}^2 e_k e_{2-k} + 2\im e_0 q_0^4 \dd{q_0}{w}.
 \end{alignat}
\end{subequations}
The general equation at $\mathcal{O}(\epsilon^{n/2})$ is unwieldly, but we shall seek only those terms that are needed to describe the $n \to \infty$ behaviour. The relevant terms are given by
\begin{multline}
 1 - \im q_0\Bigl[e_0 \phi_0 \Bigr] \frac{\phi'_{n-2}}{\phi_n} - 
 \im q_0\Bigl[e_0 \phi_1 + e_1 \phi_0 \Bigr] \frac{\phi'_{n-3}}{\phi_n} \\ - 
 \im q_0\Bigl[e_0 \phi_2 + e_1 \phi_1 + e_2 \phi_0 \Bigr]
\frac{\phi'_{n-4}}{\phi_n} - \im q_0\Bigl[e_0 \phi_0' \Bigr]
\frac{\phi_{n-2}}{\phi_n} + \ldots = 0, \label{eq:nhull_cc_Oen2}
\end{multline}
where we have divided the $\mathcal{O}(\epsilon^{n/2})$ equation by $\phi_n$ to ease the expansions to follow. As in Sec.~\ref{sec:divergence}, the merged singularity at $w = -a$ in \eqref{phi0} will cause the late terms to diverge. However, if we substitute the usual factorial-over-power ansatz into \eqref{eq:nhull_cc_Oen2} we discover terms of order $1/n^{1/2}$ are created, which cannot be matched. In fact, it can be seen through expansion of the ratios, \eg $\phi_{n-k}/\phi_n$ (\ref{sec:general}), that the divergence is described by an exponential-over-power ansatz of the form
\begin{equation} \label{eq:nhull_cc_phiansatz}
 \phi_n \sim \frac{P(w) \Gamma\left(\frac{n}{2}+\gamma\right)\exp \left[r_1(w) n^{1/2} \right]}{[\chi(w)]^{n/2+\gamma}},
\end{equation}
valid in the limit $n \to \infty$. The \emph{singulant} function, $\chi(w)$, will be solved for subject to the requirement that $\chi(-a) = 0$. Note that instead of the factorial in \eqref{eq:nhull_cc_phiansatz}, we could have equally used the ansatz
\begin{equation}
\phi_n = \frac{\tilde{P}(w)}{[\chi(w)]^{n/2 + \gamma(w)}} \exp \left[ \frac{n}{2}\log n + b(w) n + c(w) n^{1/2} + d(w) \log n\right].
\end{equation}
Had we done so, we would find $b(w) = -1/2(1 + \log 2)$, $d(w) = \gamma - 1/2$, and $\gamma(w) = \gamma$ is constant. The difference between these two ansatzes is only a numerical prefactor, with $\tilde{P}_k = \sqrt{2\pi}2^{-(\gamma - 1/2)}P$. Since we prefer to preserve the connection to the Gamma function, we continue with \eqref{eq:nhull_cc_phiansatz}.


Substituting \eqref{eq:nhull_cc_phiansatz} into (\ref{eq:nhull_cc_Oen2}) and expanding the ratios gives
\begin{multline}
 1 - \im q_0\Biggl[e_0 \phi_0 \Biggr] \Biggl\{ -\chi' + 
 \frac{\chi'r_1 + 2\chi r_1'}{n^{1/2}} + 
 \left(\frac{4\chi P' - \chi' P r_1^2 - 4\chi P r_1 r_1'}{2P}\right)
 \frac{1}{n} + \mathcal{O}\biggl(\frac{1}{n^{3/2}} \biggr) \Biggr\}
 \\ - \im q_0\Biggl[e_0 \phi_1 + e_1 \phi_0 \Biggr] 
 \Biggl\{ -\frac{\sqrt{2\chi}\chi'}{n^{1/2}} + 
 \left( \frac{3\sqrt{\chi}\chi'r_1}{\sqrt{2}} + 2\sqrt{2}\chi^{3/2} r_1'
 \right) \frac{1}{n} + \mathcal{O}\biggl(\frac{1}{n^{3/2}} \biggr)
\Biggr\} \\ - 
 \im q_0\Biggl[e_0 \phi_2 + e_1 \phi_1 + e_2 \phi_0 \Biggr] 
 \Biggl\{ -\frac{2\chi\chi'}{n} + \mathcal{O}\biggl(\frac{1}{n^{3/2}}
\biggr) \Biggr\} 
 -\im q_0\Biggl[e_0 \phi_0' \Biggr]\Biggl\{\frac{2\chi}{n} +
\mathcal{O}\biggl(\frac{1}{n^{3/2}}
\biggr) \Biggr\} = 0. \label{eq:cc_Oen2_expand}
\end{multline}
\VerbatimFootnotes
We note that these unwieldy expansions (which require higher-order contributions from Stirling's approximation to the Gamma function) can be easily derived%
\footnote{For instance, using Mathematica 10.0, we obtain the expanded ratio $\phi_{n-2}'/\phi_n$ using \\
\verb|phi[n_] = P[w] Gamma[n/2 + gamma] Exp[r1[w] Sqrt[n]]/chi[w]^(n/2 + gamma);| \\
\verb|Series[D[phi[n - 2], w]/phi[n], {n, Infinity, 1}]|}
using a computer algebra system.

Once we substitute the early orders \eqref{earlyorders} into the $\mathcal{O}(\epsilon^{n/2})$ equation \eqref{eq:cc_Oen2_expand} and take $n \to \infty$, we obtain, from the first three orders, three differential equations for the component functions of the late-order ansatz \eqref{eq:nhull_cc_phiansatz}. These are given by 
\begin{subequations}
\begin{alignat}{4}
\Oh(1)\!&:& 
\dd{\chi}{w} &= \frac{\im}{q_0^3}, \label{eq:nhull_cc13_chip} \\ 
\Oh(n^{-1/2})\!&:&
\dd{r_1}{w} &= \frac{3\im e_1}{\sqrt{2\chi} q_0^3} - \frac{\im r_1}{2\chi
q_0^3}, \label{r1p} \\
\Oh(n^{-1})\!&:& \qquad 
\frac{1}{P}\dd{P}{w} &= -\frac{4 q_0'}{q_0} + \frac{\im}{q_0^3} \left( -6 e_1^2 + 3\im e_2
+ \frac{3 e_1 r_1}{2\sqrt{2\chi}} - \frac{ r_1^2}{4\chi} \right).
\label{eq:nhull_cc13_Prat}
\end{alignat}
\end{subequations}

We substitute the value of $q_0$ in \eqref{phi0} into (\ref{eq:nhull_cc13_chip}), and solve for the singulant function subject to $\chi(-a) = 0$, giving
\begin{equation} \label{chiexplicit}
\chi(w) = \im \Bigl[ (w + a) + a\log w - a \log(-a)\Bigr] = 
a\pi + \im\Bigl[ (w+a) + a\log(w/a)\Bigr].
\end{equation}
We expect that the Stokes Phenomenon will be associated with exponentials of the form $\e^{-\chi/\ep}$ [see \eqref{2hullphiexp} and later \eqref{phiexp}], and thus we select the principal branch of the logarithm in \eqref{chiexplicit}. This ensures that $\Re(\chi) > 0$ when $w > 0$, producing exponential decay as $\epsilon \to 0$. We can then use this explicit form of $\chi$ to solve the equation for $r_1$ in \eqref{r1p}, giving
\begin{equation} \label{r1}
r_1(w) = \frac{3}{\sqrt{2\chi}} \left[\int_{-a}^w \dd{\chi}{\varphi} e_1(\varphi) \, \de{\varphi} + \text{const.}\right]
= \frac{3 \im \beta f_1}{\sqrt{2\chi}} \log (-w/a),
\end{equation}
where $e_1$ follows from \eqref{en}, and the constant of integration follows from the requirement that $r_1$ is bounded as $w \to -a$. At this point, it can be difficult to see which logarithmic and square root branch of \eqref{r1} must be used. We leave these choices ambiguous for now, but will return in the next section to clearly define the branch structure.

Turning to \eqref{eq:nhull_cc13_Prat}, we can use the preceding results for $\chi$ and $r_1$, and the expressions for $e_1$ and $e_2$ in \eqref{en} to integrate the right hand-side explicitly, giving
\begin{equation} \label{P}
P = \Delta \Bigl[q_0(w)\Bigr]^{-4 - 3\im A} \exp\Bigl[\frac{r_1^2(w)}{4}\Bigr],
\end{equation}
\noindent for some complex-valued constant $\Delta$ and where $A = 3\beta^2(2 f_1^2 - f_2)/a$.

The value of $\gamma$ in the late-order ansatz \eqref{eq:nhull_cc_phiansatz} can be derived by requiring that the behaviour of the ansatz at $n = 0$ matches the behaviour of $q_0$ in \eqref{phi0} as the singularity is approached. As $w \to -a$, we see that $\phi_n = \Oh(w+a)^{4/3 + \im A - 2\gamma}$ must match $q_0^2 = \Oh(w+a)^{-2/3}$. Solving for $\gamma$, we obtain
\begin{equation} \label{gamma_out}
	\gamma = \gamma_r + \im \gamma_c, \quad \text{where } \gamma_r = 1, \gamma_c = \frac{3\beta^2}{2a} (2 f_1^2 - f_2).
\end{equation}
To derive the above, we have made use of the limiting behaviours of $q_0$, $\chi$, $r_1$, and $P$ as the $w \to -a$. These are given by
\begin{subequations} \label{chir1P_limit}
\begin{alignat}{4}
q_0 &\sim c(w+a)^{1/3} 
\qquad &&\text{\quad where } 
c &&= (-a)^{1/3} \label{q0in} \\
\chi &\sim X(w+a)^2 
\quad &&\text{\quad where } 
X &&= -\im/2a \label{chi_in} \\
r_1 &\sim \mu_1, 
\quad &&\text{\quad where } 
\mu_1 &&= 3 \sqrt{2X} \beta f_1 \label{r1_in} \\ 
P &\sim \Delta \bigl[c (w + a)^{-\frac{1}{3}}\bigr]^{2(1 - 3\gamma)}\e^{\mu_1^2/4} 
\quad &&\text{\quad where } 
\gamma &&= 1 + \im (3\beta^2/a)(2f_1^2-f_2).
\end{alignat}
\end{subequations}

This completes our derivation of the form of the late terms \eqref{eq:nhull_cc_phiansatz} of the asymptotic expansion. In Sec.~\ref{sec:smooth}, we shall discover a connection between the exponentials switched-on due to the Stokes Phenomenon and these late terms. However, the value of $\Delta$ in \eqref{Pk} is still unknown, and is obtained by matching the inner and outer expansions.

\subsection{Step 2: Match with the solution near the singularity}

The inner region occurs near $w = -a$ when $\sqep \phi_1 \sim \phi_0$ and here, the regular asymptotic expansion \eqref{phiseries13} rearranges. Using the leading-order scaling of $\phi \sim q_0^2$, we introduce inner variables, $z$ and $\widehat{\phi}$, defined by
\begin{equation} \label{innerscale13}
w + a = (\epsilon/X)^{1/2} z \quad \text{and} \quad
\phi(w) = c^2(w+a)^{-2/3} \widehat{\phi}(z).	
\end{equation}
where $c$ and $X$ are from \eqref{q0in} and \eqref{chi_in}. We also write $q_s(w) = q_0 \widehat{q}_s$, and from \eqref{qs13}, 
\begin{equation} \label{qshat}
\widehat{q}_s = \sum_{n=0}^\infty \frac{\widehat{e}_n}{z^n}.
\end{equation}
where we have defined $\widehat{e}_n = \beta^n X^{n/2}f_n$. Using the scalings \eqref{innerscale13} in (\ref{phi_ode}), we obtain the leading-order inner equation
\begin{equation} \label{innereqn1313}
\frac{1}{2 z^2} \Bigl[ -\frac{2}{3} \widehat{q_s} \widehat{\phi}^2 + z \widehat{q_s} \widehat{\phi} \dd{\widehat{\phi}}{z} \Bigr] + \Bigl[\widehat{\phi} - \widehat{q_s}^2 \Bigr]
= 0.
\end{equation}

In order to study the leading-order solution, $\widehat{\phi}$, as it tends outwards towards the outer region, we expand
\begin{equation} \label{phiseries_in}
\widehat{\phi}(z) = \sum_{n=0}^\infty \frac{A_n}{z^n}, 	
\end{equation}
as $z \to \infty$, and substitute this series into the inner equation (\ref{innereqn1313}), giving the recurrence relation
\begin{subequations} \label{recreln1313}
\begin{gather}
A_0 = 1, \qquad
A_1 = \sum_{j=0}^1 \widehat{e}_j\widehat{e}_{1-j} \label{A0A1} \\
A_n = \sum_{j=0}^n \widehat{e}_j \widehat{e}_{n-j} + 
 \sum_{k=0}^{n-2} \widehat{e}_k \Biggl[ \sum_{j=0}^{n-2-k}\left( \frac{j + 2/3}{2} \right)  A_j A_{n-2-k-j} \Biggr] \quad \text{for $n > 1$}. \label{Ancc}
\end{gather}
\end{subequations}

A simple numerical computation assures us that $A_n$ diverges in the limit $n\to \infty$. As in the outer analysis of the previous section, a standard factorial divergence is insufficient, and an additional exponential growth in $n$ accompanies $A_n$. We posit that in the limit $n \to \infty$, the divergence of $A_n$ is captured by the two possible ansatzes
\begin{subequations} \label{An_in_both}
\begin{equation} \label{An_in}
 A_n \sim \Omcc \e^{\im \tau} \Gamma\left(\frac{n}{2} + \gamma \right)  \exp \Bigl[ \mu_1 n^{1/2} \Bigr]
\end{equation}
or
\begin{equation} \label{An_in_alt}
 A_n \sim (-1)^n \Omcc \e^{\im \tau} \Gamma\left(\frac{n}{2} + \gamma \right)  \exp \Bigl[ \mu_1 n^{1/2} \Bigr],
\end{equation}
\end{subequations}
where $\gamma$ and $\mu_1$ are real constants to be determined (and will be found to be identical to the same constants as defined in the outer analysis of Sec.~\ref{sec:outercc}). Notice also that changing the square-root branch within $\widehat{e}_n$ in \eqref{qshat} produces a factor of $(-1)^n$ in the expansion of $\widehat{q}_s$. We will proceed with the assumption that \eqref{An_in} is correct, and demonstrate how this choice affects the analysis. The remaining $\Omcc$ and $\tau$ are real constants to be determined. Once \eqref{An_in} is substituted into the recurrence relation for $A_n$ in \eqref{recreln1313}, the relevant terms are found to be
\begin{multline} \label{recAn}
 1 - \Biggl[\frac{\widehat{e}_0 A_0}{2} \Biggr] \frac{nA_{n-2}}{A_n} - 
 \Biggl[\frac{\widehat{e}_0 A_1 + \widehat{e}_1 A_0}{2} \Biggr]
\frac{nA_{n-3}}{A_n} \\ - 
 \Biggl[\frac{\widehat{e}_0 A_2 + \widehat{e}_1 A_1 + \widehat{e}_2 A_0}{2}
\Biggr] \frac{nA_{n-4}}{A_n} + 
 \Biggl[\frac{\widehat{e}_0 A_0}{3} \Biggr] \frac{A_{n-2}}{A_n} + \ldots = 0,
\end{multline}
where we have divided by $A_n$ for ease of the expansion procedure. Expanding the ratios of $A_n$ gives
\begin{multline}
 1 - \Biggl[\frac{\widehat{e}_0 A_0}{2} \Biggr] 
 \Biggl\{ 2 - \frac{2\mu_1}{n^{1/2}} + \frac{4-4\gamma + \mu_1^2}{n} +
\mathcal{O}\biggl(\frac{1}{n^{3/2}} \biggr) \Biggr\} \\
 - \Biggl[\frac{\widehat{e}_0 A_1 + \widehat{e}_1 A_0}{2} \Biggr] 
 \Biggl\{\frac{2\sqrt{2}}{n^{1/2}} - \frac{3\sqrt{2}\mu_1}{n} +
\mathcal{O}\biggl(\frac{1}{n^{3/2}} \biggr) \Biggr\}  \\ - 
 \Biggl[\frac{\widehat{e}_0 A_2 + \widehat{e}_1 A_1 + \widehat{e}_2 A_0}{2}
\Biggr] 
  \Biggl\{ \frac{4}{n} + \mathcal{O}\biggl(\frac{1}{n^{3/2}}
\biggr) \Biggr\} 
  + \Biggl[\frac{\widehat{e}_0 A_0}{3} \Biggr]\Biggl\{ \frac{2}{n} +
\mathcal{O}\biggl(\frac{1}{n^{3/2}} \biggr) \Biggr\}  = 0. \label{eq:nhull_recttmp} 
\end{multline}

\noindent Setting $\widehat{e_0} = A_0 = 1$, then (\ref{eq:nhull_recttmp}) is satisfied identically at $\mathcal{O}(1)$. Afterwards, the $\mathcal{O}(1/n^{1/2})$ and $\mathcal{O}(1/n)$ terms give
\begin{subequations}
\begin{gather}
 \mu_1 = 3\sqrt{2} \widehat{e_1} = 
 3\sqrt{2 X} \beta f_1,  \label{eq:nhull_cc13_mu1} \\
 \gamma = 1 - 6 \widehat{e_1}^2 + 3 \widehat{e_2} = 
 1 - 3\beta^2 X \left[ 2f_1^2 - f_2\right]. \label{gamtemp}
\end{gather}
\end{subequations}
where we have set $\widehat{e_1} = \beta \sqrt{X} f_1$ and
$\widehat{e_2} = \beta^2 X f_2$. In fact, this independently verifies the inner limit of $r_1(w) \to \mu_1$ obtained in the outer analysis in \eqref{r1_in}, as well as the $\gamma$ constant in \eqref{gamma_out}. 

Now that the exponential growth of the divergent ansatz of $A_n$ in \eqref{An_in} has been fully determined, the values of $\Omcc$ and $\tau$ can be computed by numerically solving the full nonlinear recurrence relation \eqref{recreln1313} and using 
\begin{equation} \label{Hn}
	H_n \equiv \frac{A_n}{\Gamma(\frac{n}{2} + \gamma) e^{\mu_1 n^{1/2}}} \to \Omcc \e^{\im \tau}
\end{equation}
in the limit $n\to\infty$. 

\begin{figure}[htb] \centering
\subfigure [$\sigma_1 = 3/24$, $\sigma_2 = 5/24$] {
\includegraphics{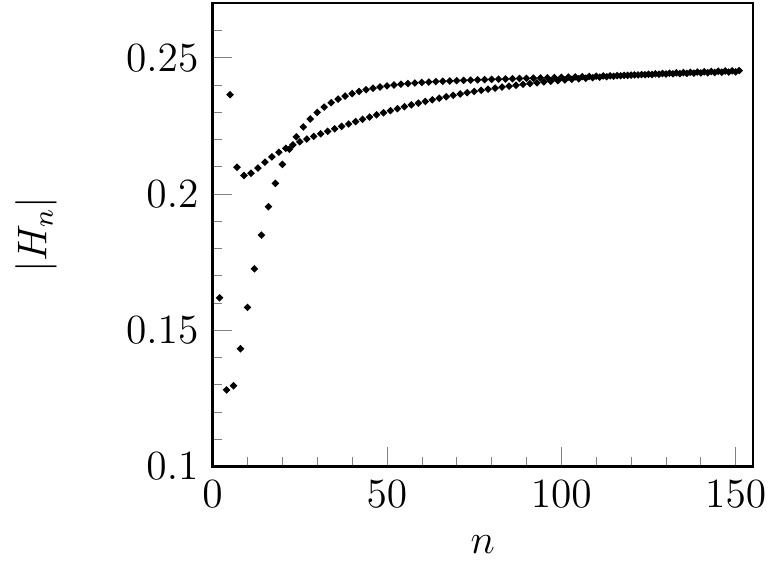}
\includegraphics{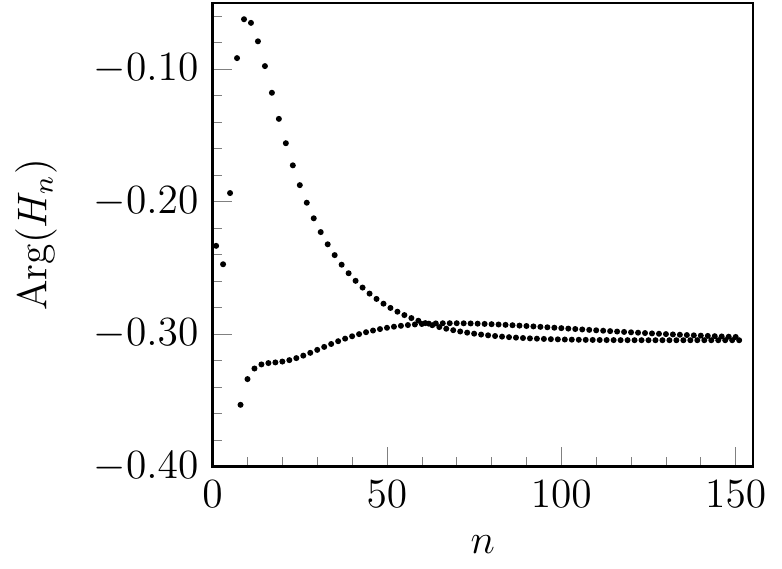}
} \\
\subfigure [$\sigma_1 = 6/24$, $\sigma_2 = 2/24$] {
\includegraphics{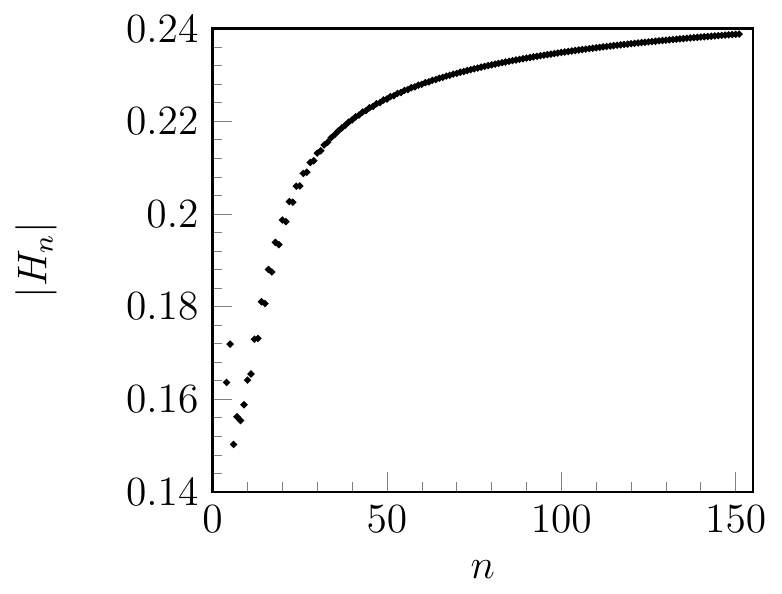}
\includegraphics{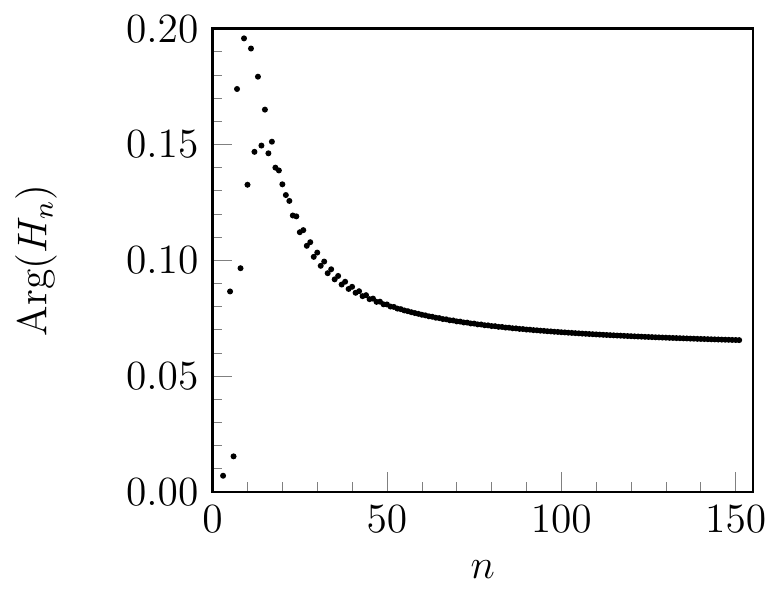}
} 
\caption{Magnitude and phase of $H_n = A_n/[\Gamma(n/2 + \gamma) \e^{\mu_1 n^{1/2}}]$, representing the solutions of the recurrence relation \eqref{recreln1313} once the exponential growth has been scaled out. Two different $\sigma_1$ and $\sigma_2$ pairs are used, and $a = 1$, $\beta = 1$ for both. \label{recplot}}
\end{figure}

In Figure \ref{recplot}, we plot $H_n$ for the two cases (a) $\sigma_1 = 3/24$, $\sigma_2 = 5/24$ and (b) $\sigma_1 = 6/24$, $\sigma_2 = 2/24$. The remaining parameters are set to $a = 1$ and $\beta = 1$. The convergence towards the constant values in each of the graphs is algebraic in $n$, and confirms that the exponential growth was correctly predicted. The values of $|H_n|$ and $\text{Arg}(H_n)$, can be seen to alternate between two branches, and this is effectively a consequence of the $\sqep$ series in $q_s$ and our choice of $\ell/m = 1/2$. For other $\ell/m$ values in \eqref{eq:lmsig}, it can be expected that $m$ branches of the recurrence relation would be observed (see \ref{sec:general} for discussion of the methodology for more general singularities). 

Now had we used \eqref{An_in_alt} with the extra factor of $(-1)^n$, we would note that the value of $\gamma$ in \eqref{gamtemp} remains the same, but the value of $\mu_1$ in \eqref{eq:nhull_cc13_mu1} changes to $-3\sqrt{2}\widehat{e}_1$. Both forms of \eqref{An_in_both} are possible, but one contribution will exponentially dominate the other. The proper choice of ansatz can be determined based on the fact we need $\Re(\mu_1) > 0$ in order to describe the divergent $A_n$ values in \eqref{Hn}. Since $f_1 = \sigma_2 - \sigma_1$ in \eqref{eq:nhull_cc_fn} and $\beta > 0$, we define 
\begin{equation} \label{mu1_branchpick}
\mu_1 = 3\beta |\sigma_2 - \sigma_1| \sqrt{2|X|} \, \e^{-\pi \im/4}, 
\end{equation}
which ensures that $\Re(\mu_1) > 0$ regardless of the sign of $(\sigma_2 - \sigma_1)$. Indeed, without taking care to choose the correct ansatz, the graphs of $|H_n|$ in Figure \ref{recplot} would display exponential growth (due to underpredicting the divergence) or the points in the graphs of $\text{Arg}(H_n)$ may alternate by $\pm \pi$. 

\subsection{Step 3: Matching inner and outer expansions}

The matching between the outer expansion \eqref{phiseries13} and inner expansion \eqref{phiseries_in} will be imposed through Van Dyke's rule \cite{vandyke_1975}. The principle states that 
\begin{multline}
\text{$m$-term inner expansion of ($n$-term outer expansion)} \\ = 
\text{$n$-term outer expansion of ($m$-term inner expansion)},
\end{multline}
or in our shorthand, $\text{($m$.t.i)($n$.t.o)} = \text{($n$.t.o)($m$.t.i)}$. The $n^\text{th}$ term of the outer expansion given by the ansatz \eqref{eq:nhull_cc_phiansatz}, written in inner coordinates, and keeping the first term yields
\begin{equation}
\phi \xrightarrow[]{\text{($n$.t.o)}} \epsilon^\frac{n}{2} \phi_n \sim
 \frac{P
 \Gamma(\frac{n}{2}+\gamma)\e^{r_1 n^{1/2}}}{\chi^{\frac{n}{2}+\gamma}} 
 \xrightarrow[]{\text{($1$.t.i)}} \frac{\Delta c^{2(1-3\gamma)}W^{-2\sigma} \e^{\mu_1^2/4}
\Gamma(\frac{n}{2}+\gamma)\e^{\mu_1 n^{1/2}}}{X^\gamma z^n},
\end{equation}
 where we have used the inner limits of $\chi$, $r_1$, and $P$ in \eqref{chir1P_limit}. This expression provides a match with the first term of the inner approximation ($1$.t.i) given by \eqref{phiseries_in}, written in outer variables and re-expanded to $n$ terms ($n$.t.o), 
\begin{equation}
\phi = c^2 W^{-2\sigma} \widehat{\phi} \xrightarrow[]{\text{($1$.t.i)}} 
 c^2 W^{-2\sigma} \sum_{n=0}^\infty \frac{A_n}{z^n} 
 \xrightarrow[]{\text{($n$.t.o)}} 
c^2 W^{-2\sigma} \frac{[\Omcc \e^{\im \tau}] \Gamma(\frac{n}{2} + \gamma)e^{\mu_1 n^{1/2}}}{z^n}.	
\end{equation}
 Matching then allows us to obtain the value of outer prefactor, $\Delta$, as a function of the inner pre-factors $\Omcc \e^{\im \tau}$, and thus $\Delta = c^{6\gamma} X^\gamma \e^{-\mu_1^2/4} \left[ \Omcc \e^{\im\tau} \right]$. In summary, $P$ is given by 
 \begin{equation} \label{Pagain}
 	P(w) = (c^{6} X)^\gamma \bigl[\Omcc \e^{\im \tau}\bigr] \bigl[q_0(w)\bigr]^{2(1 - 3\gamma)} \exp\bigl[\tfrac{1}{4}(r_1^2(w) - \mu_1^2)\bigr].
 \end{equation}
 We have thus obtained all the components of the late-order ansatz for $\phi_n$ in \eqref{eq:nhull_cc_phiansatz}. It remains to relate these late terms to the exponential switching.

\subsection{Step 4: Optimal truncation and Stokes line smoothing} \label{sec:smooth}

We now truncate the asymptotic expansion of $\phi$ in \eqref{phiseries13} at $n = \N$, 
\begin{equation}
\phi = \sum_{n=0}^{\mathcal{N}-1} \epsilon^{n/2} \phi_n + R_\mathcal{N},	
\end{equation}
and substitute this equation into the differential equation for $\phi$ in \eqref{phi_ode}. This yields an equation of the form 
\begin{subequations} \label{Rn_eqn}
\begin{equation} \label{eq:nhull_cc13_remeqn}
 \L(R_\mathcal{N}; \epsilon) \sim \epsilon^{\mathcal{N}/2} \phi_\mathcal{N},
\end{equation}

\noindent with $\L$ given by
\begin{multline} \label{eq:nhull_cc_Lop2}
 \L(R_\mathcal{N}; \epsilon) \equiv R_\mathcal{N} - \im q_0\Bigl[e_0 \phi_0 \Bigr] \epsilon
R'_\mathcal{N} - 
 \im q_0\Bigl[e_0 \phi_1 + e_1 \phi_0 \Bigr] \epsilon^{3/2} R'_\mathcal{N} \\ - 
 \im q_0\Bigl[e_0 \phi_2 + e_1 \phi_1 + e_2 \phi_0 \Bigr] \epsilon^2 R'_\mathcal{N}
 -\im q_0\Bigl[e_0 \phi_0' \Bigr] \epsilon R_\mathcal{N} + \ldots
\end{multline}
\end{subequations}

Our goal in this section is to choose $n = \N$ optimally, so that the remainder is exponentially small. In the limit that $\epsilon \to 0$, $\N \to \infty$, and the inhomogeneous differential equation for the remainder \eqref{eq:nhull_cc13_remeqn} will be forced by the exponential-over-power ansatz of \eqref{eq:nhull_cc_phiansatz}. We shall find that as the Stokes line, $\Im(\chi) = 0$ and $\Re(\chi) \geq 0$, is crossed the exponentially small remainder switches on. The distinction of our work here, compared to the previous studies of, \eg Chapman \etal \cite{chapman_1998}, is that the new exponential-over-power ansatz modifies the connection between $\phi_\N$ and $R_\N$. Effectively, the Stokes line is shifted in location due to the multiple singlarities. 

To begin, we express the solution of the homogeneous equation, $\L = 0$, in \eqref{Rn_eqn}, by setting
\begin{equation} \label{Rn_homo}
 R_\mathcal{N} = \overline{P}(w) \e^{F(w)} = 
\overline{P}(w) \exp\Biggl[\frac{F_0(w)}{\epsilon} +
\frac{F_1(w)}{\sqep}\Biggr],
\end{equation}

\noindent into (\ref{eq:nhull_cc_Lop2}) and dividing by $R_\mathcal{N}$, giving
\begin{multline}
 1 - \im q_0\Biggl[e_0 \phi_0 \Biggr] \Biggl[F_0' + \ep^{\frac{1}{2}} F_1' + \epsilon
\frac{\overline{P}'}{\overline{P}}\Biggr] - 
 \im q_0\Biggl[e_0 \phi_1 + e_1 \phi_0 \Biggr] \Biggl[ \sqep F_0' +
\epsilon F_1' + \mathcal{O}(\epsilon^{\frac{3}{2}}) \Biggr]
  \\ - 
 \im q_0\Biggl[e_0 \phi_2 + e_1 \phi_1 + e_2 \phi_0 \Biggr] 
 \Biggl[ \epsilon F_0' + \mathcal{O}(\epsilon^{\frac{3}{2}}) \Biggr]
 -\im q_0\Biggl[e_0 \phi_0' \Biggr] \epsilon = 0.
\end{multline}

\noindent Solving at each order yields
\begin{subequations}
\begin{alignat}{3}
\Oh(1)\!&:& \qquad
\dd{F_0}{w} &= -\frac{\im}{q_0^3} = - \chi' \label{dF0}, \\
\Oh(\sqep)\!&:& \qquad 
\dd{F_1}{w} &= \frac{3\im e_1}{q_0^3} = 3 e_1 \chi', \label{F1p} \\
\Oh(\ep)\!&:& \qquad 
\frac{1}{\overline{P}}\dd{\overline{P}}{w} &= 
- \frac{4q_0'}{q_0} + \frac{3\im e_2}{q_0^3} - \frac{6\im e_1^2}{q_0^3}. \label{eq:nhull_cc_Rncomp_last}
\end{alignat}
\end{subequations}

\noindent These expressions for $F_0 = -\chi$, $F_1$ and $\overline{P}$ are related to the functions $\chi$, $r_1$, and $P$ from the late-order ansatz (\ref{eq:nhull_cc_phiansatz}). First, we compare the equations for $r_1'$ \eqref{r1p} and $F_1'$ \eqref{F1p} to conclude that $\dd{}{w} \left[ \sqrt{2\chi} r_1 - F_1\right] = 0$, and thus
\begin{equation} \label{F1_r1}
F_1(w) = \sqrt{2\chi} r_1,
\end{equation}
where we have set the constant of integration so that $F_1 = 0$ at the singularity, $w = -a$, where $\chi = 0$ and $r_1 = \mu$. Also, notice from \eqref{eq:nhull_cc13_Prat} and \eqref{eq:nhull_cc_Rncomp_last} that $P$ and $\overline{P}$ are related through 
\begin{equation} \label{Prelate}
\overline{P}(w) = P(w) \exp\left[-\frac{r_1^2(w)}{4}\right].
\end{equation}

In order to solve the inhomogeneous equation, we multiply \eqref{Rn_homo} by the Stokes Smoothing parameter, so that $R_\mathcal{N} = \St\overline{P} e^F$, and substitute into (\ref{Rn_eqn}), giving
\begin{equation}
-\epsilon \im q_0^3 \dd{\St}{w} \overline{P} \exp\left[\frac{F_0}{\ep} + \frac{F_1}{\sqep}\right] 
\sim \epsilon^{\frac{\mathcal{N}}{2}}
\frac{P \Gamma\left(\frac{\mathcal{N}}{2}+\gamma\right)
\e^{r_1 n^{\frac{1}{2}}}
}{\chi^{\frac{\mathcal{N}}{2}+\gamma}}.	
\end{equation}
For this expression, we write the derivative in terms of $\chi$ so that $\de\St/\de{w} = \chi' \de\St/\de{\chi}$ giving
\begin{equation} \label{eq:nhull_ssmooth_tmp}
\epsilon \dd{\St}{\chi} \Bigl[ \overline{P} e^{F_1/\sqep}
\Bigr] \sim
\Bigl[Pe^{r_1 \mathcal{N}^{\frac{1}{2}}} \Bigr] 
\frac{ \epsilon^{\frac{\mathcal{N}}{2}} \Gamma\left(\frac{\mathcal{N}}{2}+\gamma\right)
e^{\frac{\chi}{\epsilon}}}
{\chi^{\frac{\mathcal{N}}{2}+\gamma}}.
\end{equation}
The optimal truncation point is found where adjacent terms in the asymptotic approximation are of the same size, and thus $\epsilon^{(\mathcal{N}+1)/2}
|\phi_{\mathcal{N}+1}| \sim \epsilon^{\mathcal{N}/2}
|\phi_\mathcal{N}|$. Using the ansatz \eqref{eq:nhull_cc_phiansatz} and writing 
\begin{equation} \label{chi_rnu}
\chi = r \e^{\im \nu}	
\end{equation}
we find that the optimal truncation point is found where
\begin{equation} \label{eq:nhull_cc13_Nval}
 \mathcal{N} = \frac{2r}{\epsilon} + 2\rho,
\end{equation}
with $\rho \in (0, 1/2)$ as $\epsilon \to 0$. Our goal now is to expand the equation for the Stokes smoothing factor, $\mathcal{S}$, in \eqref{eq:nhull_ssmooth_tmp} in the limit $\epsilon \to 0$, and examine its rate of change as the Stokes line is crossed (a fixed $r$ and varying $\nu$ in \eqref{chi_rnu}, which corresponds to the Stokes line $\Re(\chi) \geq 0$). We first use Stirling's formula to write
\begin{equation} \label{Gamma_expand}
\Gamma(\mathcal{N}/2 + \gamma) \sim 
\sqrt{2\pi} \left(\frac{r}{\epsilon}\right)^{\mathcal{N}/2 + \gamma - 1/2} \e^{-r/\epsilon} \left(1 + O(\epsilon)\right).	
\end{equation}
It next follows from the relation between $r_1$ and $F_1$ in \eqref{F1_r1} that
\begin{equation} \label{r1_expand}
r_1 \mathcal{N}^{1/2} - \frac{F_1}{\sqep} 
= r_1 \mathcal{N}^{1/2} - \frac{\sqrt{2\chi}r_1}{\epsilon}= r_1 \sqrt{\frac{2r}{\epsilon}} \left(1 - e^{i\nu/2} + \frac{1}{2} \frac{\epsilon}{r} \rho + O(\epsilon^2)\right).
\end{equation}
We can now exchange the differentiation in $\chi$ to differentiation across the Stokes line, in $\nu$, using  
$\de{\St}/\de{\chi} = -(\im \e^{-\im\nu}/r) \, \de{\St}/\de{\nu}$. Using the substitution \eqref{Gamma_expand} for $\Gamma$ and \eqref{r1_expand} for $r_1$ in the equation for $\St$ in \eqref{eq:nhull_ssmooth_tmp} now yields
\begin{equation} \label{dSt_dnu}
\dd{\St}{\nu} \sim \im \left(\frac{P}{\overline{P}}\right) \frac{\sqrt{2\pi r}}{\epsilon^{\gamma + 1/2}} \e^G,
\end{equation}
where we have defined $G$ according to
\begin{equation} \label{G_expand}
G = -\frac{r}{\epsilon}\left(1 - \e^{\im\nu} + \im\nu\right) + 
r_1 \sqrt{\frac{2r}{\epsilon}} \left(1 - \e^{\im\nu/2} + \frac{1}{2} \frac{\epsilon}{r} \rho + O(\epsilon^2)\right) + \im\nu (1 - \rho - \gamma)	
\end{equation}
The first group of bracketed terms in the above expression for $G$ indicate that the Stokes smoothing constant, $\St$, in \eqref{dSt_dnu} is exponentially small unless $\nu$ is also small. We rescale $\nu = \sqep \overline{\nu}$ and note that $G = -r\overline{\nu}^2/2 - \im r_1\sqrt{r/2} + \Oh(\sqep)$. Substitution of this expression into the equation for $\St$ in \eqref{dSt_dnu} yields after simplification
\begin{equation}
\dd{\St}{\overline{\nu}} \sim 
\im \left(\frac{P}{\overline{P}}\right) \frac{\sqrt{2\pi r}}{\epsilon^{\gamma}} \exp\left[ -\left(\sqrt{\frac{r}{2}}\overline{\nu} + \im \frac{r_1}{2}\right)^2\right] \e^{-r_1^2/4}.
\end{equation}
Finally, we may use the relationship between $P$ and $\overline{P}$ in \eqref{Prelate} to conclude that 
\begin{equation}
\dd{\St}{\overline{\nu}} \sim \im \frac{\sqrt{2\pi r}}{\epsilon^{\gamma}} \exp\left[ -\left(\sqrt{\frac{r}{2}}\overline{\nu} + \im \frac{r_1}{2}\right)^2\right].
\end{equation}

It remains now to only substitute $\overline{\theta} = \overline{\nu} + \im r_1/\sqrt{2r}$ in order to properly centre the integration over the Stokes line. Because $r_1$ is locally constant near the Stokes line, and because $|\chi| = r$ is fixed (by the optimal truncation) then we must have $\de{\overline{\theta}}/\de{\overline{\nu}} = 1 + \Oh(\sqep)$. Thus, 
\begin{equation} \label{dSt_final}
  \dd{\St}{\overline{\theta}} \sim \im \frac{\sqrt{2\pi r}}{\epsilon^{\gamma}} \e^{-r \overline{\theta}^2/2}.
\end{equation}
which is precisely the same expression for the Stokes smoothing that occurs for the case of well-separated singularities \cite{chapman_2006}. We integrate \eqref{dSt_final} from $\overline{\theta} = \infty$ to $\overline{\theta} = -\infty$, which corresponds to crossing the Stokes line as $w$ increases [due to the geometry of the Stokes line in Figure \ref{fig:stokes} coupled with parameterization of $\chi$ in \eqref{chi_rnu}]. We thus obtain $\St \sim -2\pi\im/\epsilon^\gamma$. To conclude, the remainder switched-on across Stokes lines is
\begin{equation} \label{phiexp}
\phicc \equiv \Bigl[ R_\mathcal{N} \Bigr]_-^+ \sim -\frac{ 2\pi \im}{\epsilon^\gamma} \left[P(w) \e^{-r_1^2(w)/4}\right] \exp \left[-\frac{\chi}{\epsilon} + \frac{F_1}{\sqep}\right]	
\end{equation}
moving in the direction of increasing $w$, and where we have used the expression for $F_0$ in \eqref{dF0}, and the replacement of $\overline{P}$ by $P$ in \eqref{Prelate}. 

We wish to obtain an expression for $|\phicc|$ when $w$ is real and positive. Recall that the branch of $\chi$ was chosen in \eqref{chiexplicit} so that $\Re(\chi) = a > 0$ for $w > 0$, which indeed corresponds to exponential decay in \eqref{phiexp}. Similarly, we select the branch of $F_1 = \sqrt{2\chi} r_1$ in \eqref{F1_r1} so that, it too, will be assured to possess a positive real value on $w > 0$. This yields
\begin{equation}
F_1(w) = - 3\pi \beta |\sigma_2 - \sigma_1| + \im \bigl[ 3\beta (\sigma_2 - \sigma_1) \log(w/a) \bigr],
\end{equation}
where we have used $r_1$ in \eqref{r1} and $f_1 = \sigma_2 - \sigma_1$. 

Next, with $P$ given in \eqref{Pagain} and $\gamma = \gamma_r + \im \gamma_c = 1 + \im \frac{3\beta^2}{2a}(2 f_1^2 - f_2)$ from \eqref{gamma_out}, we have
\begin{equation} \label{Pe_mag}
|P(w) \e^{-r_1^2(w)/4}| = |c^{6} X| \e^{-\frac{9\beta^2}{4a}(2f_1^2 - f_2)}\Omcc |q_0(w)|^{2(1 - 3\gamma)},
\end{equation}
since $\mu_1^2$ appearing in $P$ purely imaginary by \eqref{r1_in}. Note that in the above formula, we have also used the fact %
that $\gamma_c \text{Arg}(c^6 X) = \gamma_c \text{Arg}(-\im) = \gamma_c(3\pi/2)$. Using $|c^6 X| = a/2$ from \eqref{q0in} and \eqref{chi_in}, and combining \eqref{phiexp}--\eqref{Pe_mag}, we conclude that the amplitude of the exponentially small waves is 
\begin{equation} \label{phiexp_abs}
\left\lvert \phicc \right\rvert \sim \left[\frac{\pi a \Omcc}{\epsilon q_0^4(w)}\right] \exp\left[-\frac{9\beta^2}{4a}(2f_1^2 - f_2)\right]  \exp \left[-\frac{a\pi}{\epsilon} - \frac{3\pi \beta |\sigma_2 - \sigma_1|}{\sqep}\right].
\end{equation}

\section{Numerical verification of $\sigma_1 = 1/6$ and $\sigma_2 = 1/6$}
\label{sec:nhull_cc1616}

\noindent We will now apply the formulae developed over the
course of the previous sections to the case where $q_s$ contains singularity powers of $\sigma_1 = 1/6$ and $\sigma_2 = 1/6$. In addition to comparing the asymptotic predictions to numerical computations, the most important behaviour to verify is that, in the limit $\beta \to 0$, the close-singularity approximation tends to the one-singularity approximation, and as $\beta \to \infty$, the close-singularity approximation tends to the two-singularity approximation. Our choice of the symmetric case of $\sigma_1 = \sigma_2 = 1/6$ helps to simplify some of the needed computations. In particular, $f_1 = 0$ and $f_2 = 1/6$.

To begin, we recall that there are two forms of the $q_s$ function in the differential equation \eqref{phi_ode} of interest: the two-singularity \eqref{qs_2hull} and single-singularity \eqref{qs_one} versions given by
\begin{equation}
q_s = \frac{w^{1/3}}{(w + a_1)^{1/6} (w + a_2)^{1/6}} \quad \text{and} \quad 
q_s = \left(\frac{w}{w+a}\right)^{1/3}. 
\end{equation}
The well-separated case, where $a_2 - a_1 = \Oh(1)$, was reviewed in Sec.~\ref{sec:wellsep}. The wave amplitudes for both well-separated and single-singularity cases are given by \eqref{phiexp_two}, and depend on a crucial pre-factor, $\Omega$, calculated through a numerical solution of a nonlinear recurrence relation. Both the well-separated and single-singularity variants use the same $\Omega$, whereas the closely-separated analysis in Sec.~\ref{sec:cc13} requires an $\Omcc$ that depends on $\beta$. We shall write
\begin{equation}
	\Omega(\sigma_k) \quad \text{and} \quad \Omcc(\sigma_1, \sigma_2; \beta),
\end{equation}
for the two versions, where $\sigma_k$ is the local singularity power which corresponds to the generated exponential.

For the close-singularity case, with $a_1 = a + \sqep\beta$ and $a_2 = a - \sqep\beta$, the wave amplitudes are given by \eqref{phiexp_abs}, and in the limit $w \to \infty$, $q_0 \to 1$, so we have
\begin{equation} \label{eq:nhull_cc1616_qcc}
 \left\lvert  \phicc  \right\rvert 
 = a \exp\left[\frac{3\pi
\beta^2}{8a}\right] \Omcc(\tfrac{1}{3}; \beta) \left[ \frac{\pi}{\epsilon} \e^{-\frac{\pi a}{\epsilon}} \right].
\end{equation}

\begin{figure}[htb] \centering
\includegraphics{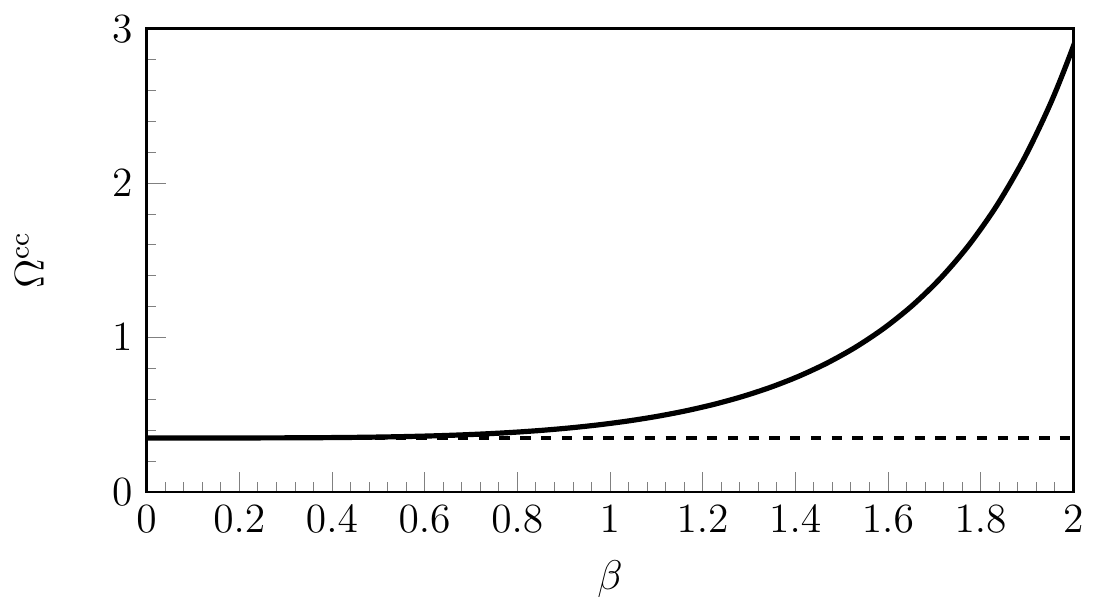}
\caption{The solid curve is the pre-factor $\Omcc(1/6,1/6; \beta)$ for the close-singularity problem, with $a = 0.5$. As $\beta \to 0$ and the two singularities coalesce, the pre-factor approaches the same value as for the one-singularity case, with $\Omega(1/3) \approx 0.351$ (dashed). \label{fig:nhull_beta0}}
\end{figure}

For the single-singularity case, we have $q_s \sim c(w+a)^{-1/3}$ near the singularity with $c = (-a)^{1/3}$. We thus apply the amplitude approximation \eqref{phiexp_two} with $\sigma_k = 1/3$, $c_k = c$, and $\gamma_k = 1$, which follows from \eqref{gammak}. Since the outer analysis of the close-singularity analysis involves a derivation of the singulant, $\chi$, corresponding to a single merged singularity, then the $\chi = \chi_k$ is given by \eqref{chiexplicit}, and $\Re(\chi) = a$ along the positive real axis. This yields the amplitude estimate
\begin{equation} 
 \left\lvert \phisingle \right\rvert \sim a \Omega(\tfrac{1}{3}) \left[\frac{\pi}{\epsilon} \e^{-\frac{\pi a}{\epsilon}} \right].
\end{equation}
 Thus, in the limit $\beta \to 0$, we have $\Omcc_1\left(\tfrac{1}{3}; \beta\right) \to \Omega\left(\frac{1}{3}\right)$ as expected. This limiting behaviour is shown in Figure \ref{fig:nhull_beta0} for the case of $a = 0.5$ ship where we see that indeed, the close-singularity approximation tends to the one-singularity approximation as $\beta \to 0$.

We now turn to the well-separated two-singularity approximation. Here, there are two singulant functions (and hence two exponentially small waves) given by 
\begin{equation}
\chi_1(w) = \int_{-a_1}^w \frac{\im}{q_s^3(\varphi)} \, \de{\varphi} \quad \text{and} \quad 
\chi_2(w) = \int_{-a_1}^w \frac{\im}{q_s^3(\varphi)} \, \de{\varphi}.
\end{equation}
However, it can be verified that if $\sigma_2 = 1/6$ and $\sigma_1$ approaches $1/6$ from above, then the previous Stokes line from $w = -a_1$ has flattened onto the real $w$ axis and now lies between $-a_1 \leq w \leq -a_2$. This follows by virtue of $\chi' = \im/q_s^3 > 0$ in this region. Thus, we conclude that waves are entirely generated by the singularity at $w = -a_2$. For $w > 0$, we have
\begin{equation} \label{rechi2}
\Re(\chi_2) = \Re \left(\int_{-a_2}^{-a_1} + \int_{-a_1}^w \right) \frac{\im}{q_s^3(\varphi)} \, \de{\varphi} = \pi a - \pi a\sqrt{1 - \epsilon (\beta/a)^2},
\end{equation} 
where the first term on the right hand-side follows from $\Re(\chi_1) = \pi (a_1 + a_2)/2 = \pi a$ from a residue contribution at infinity. Thus, the wave from the small perturbation about the point $w = -a$ has produced an additional exponential factor of $\e^{\pi\beta^2/2a}$. 

\begin{figure}[htb] \centering
\includegraphics{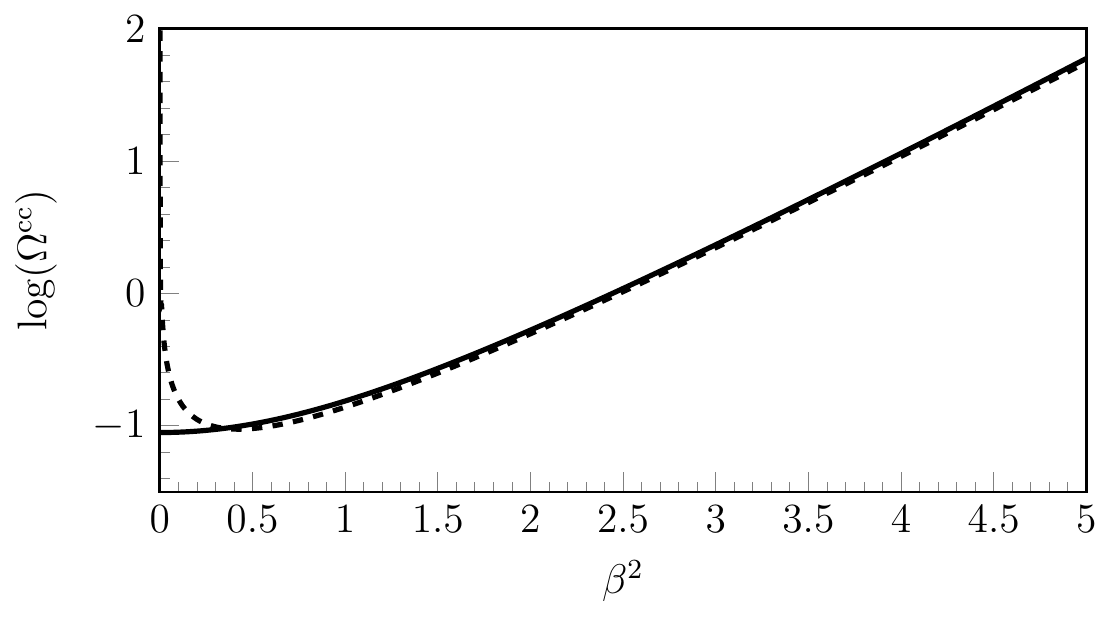}
\caption{The solid curve is the pre-factor $\Omcc(1/6,1/6; \beta)$ for close-singularity approximation, with $a = 0.5$. As $\beta \to \infty$ and the singularities separate, the pre-factor approaches the dashed curve with $\Omega(1/6) \times 2[a/(9\beta^2)]^{1/3} \exp[\pi \beta^2/3a]$, corresponding to the well-separated result. \label{fig:nhull_betainf}}
\end{figure}

With $q_s \sim c_2(w + a_2)^{\sigma_2}$ near the singularity, we obtain values of 
\begin{equation}
 \sigma_2 = 1/6, \qquad c_2 = \frac{a^{\frac{1}{3}} \e^{\pi \im/3}}{(2\sqep\beta)^{\frac{1}{6}}}, \qquad
 \gamma_2 = 2/3
\end{equation}
where $\gamma_2$ follows from \eqref{gammak}. Using \eqref{rechi2} in \eqref{phiexp_two}, we see that the amplitude of the exponentially small waves from the well-separated analysis yields
\begin{equation} \label{eq:nhull_cc1616_qtc}
\left\lvert \phi_\text{exp, $2$} \right\rvert
\sim \frac{2a^{\frac{4}{3}} \e^{\frac{\pi\beta^2}{2a}}}{(3\beta)^{\frac{2}{3}}} 
\Omega(\tfrac{1}{6}) \left[ \frac{\pi}{\epsilon} \e^{-\frac{\pi a}{\epsilon}} \right].
\end{equation}
Now comparing with the close-singularity approximation in (\ref{eq:nhull_cc1616_qcc}), we see that in order to match with (\ref{eq:nhull_cc1616_qtc}) in the limit $\beta \to \infty$, we require
\begin{equation} \label{eq:nhull_cc1616_betainf}
 \Omcc\left(\tfrac{1}{6},\tfrac{1}{6}; \beta\right) \sim
2\left(\frac{a}{9\beta^2}\right)^{1/3} \e^{\frac{\pi \beta^2}{8a}} \Omega\left(\tfrac{1}{6}\right) .
\end{equation}
In Figure \ref{fig:nhull_betainf}, we plot the natural logarithms of
the left and right-hand sides of (\ref{eq:nhull_cc1616_betainf}) as a
function of $\beta^2$ and indeed, the convergence between the two values is very fast.

\begin{figure}[htb] \centering
\includegraphics{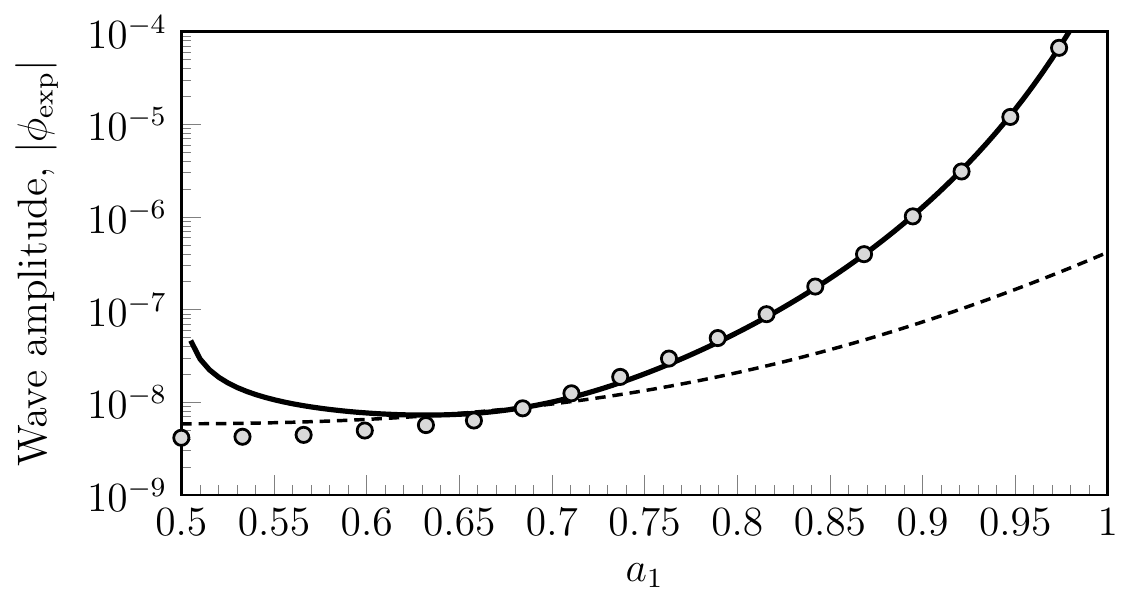}
\caption{The numerical solution (nodes) for the $\sigma_1 = \sigma_2 = 1/6$ forcing versus the asymptotic approximations for $\epsilon = 0.075$ and $a_1 + a_2 = 1$. The two-singularity approximation (solid line) is singular when $a_1, a_2 \approx 0.5$ near the left. Here, the solution is well approximated by the close-singularity approximation (dashed). \label{fig:nhull_cc1616}}
\end{figure}

The final result is shown in Figure \ref{fig:nhull_cc1616}; here, we compare the numerical and asymptotic approximations for far field wave amplitudes of the $\sigma_1 = \sigma_2 = 1/6$ forcing. As expected, the two-singularity approximation of (\ref{eq:nhull_cc1616_qtc}) for the wave amplitude is a fine approximation, but only until the two singularities begin merging near $a_1 = a_2 = 0.5$. At this point, the close-singularity approximation of (\ref{eq:nhull_cc1616_qcc}) provides a much better match to the numerical results.

\section{Discussion} \label{sec:discuss}

In this paper, we proposed a model singularly perturbed differential equation that was inspired by studies on water waves and ship hydrodynamics \cite{trinh_2014}. For this toy model, the divergence of the late terms of the associated asymptotic expansion was not described through the common factorial-over-power ansatz (Dingle \cite{dingle_1973}), but instead through a more general exponential-over-power ansatz. By applying methods in exponential asymptotics, optimally truncating the series, and smoothing the Stokes line, we were able to recover the exponentially small waves switched on. 

Now, in regards to the original motivation of water waves and hydrodynamics, the limit of coalescing singularities represents a rather niche area of practical interest; indeed, one would argue that the effort to derive the final exponential $\phicc$ in \eqref{phiexp_abs} far exceeds the effort to numerically solve the differential equation directly! However, the more important lesson from this work is in relation to the general class of nonlinear problems, $\Nop(z, y; \epsilon) = 0$, introduced in \eqref{Nop}. For such problems where the solution of the differential equation is forced by two interleaving asymptotic expansions [\eg the expansion resulting from the $\epsilon$ forcing in \eqref{phi_ode}, and the expansion resulting from the $\epsilon^{\ell/m}$ in \eqref{eq:nhull_cc_qs}], we must expect the use of more generalized exponential-over-power divergence. Our methodology in this paper highlights several interesting aspects of the adjusted exponential asymptotics theory, including the merging of multiple Stokes lines in an outer region, combined with a thickening of such lines during the optimal truncation procedure.

We may also ask the question of why such exponential-over-power divergence has not been necessary for similar problems of coalescing singularities in the exponential asymptotics literature. This situation of interleaving terms also arises in the Saffman-Taylor viscous fingering problem (see for example, \cite{combescot_1988} and \cite{chapman_1999}). If we follow the same ideas as presented there, then we would expect our expansion for $\phi$ to split into $m$ sub-expansions: 
\begin{multline} \label{eq:nhull_phiasym_split}
 \phi = \biggl[ \phi_0 + \epsilon^{\frac{1}{m}} \phi_1 + \epsilon^{\frac{2}{m}} \phi_2 + \ldots \biggr] 
 + \sum_{n=n^*}^\infty \epsilon^{n} \phi_{mn} + 
 \sum_{n=n^*}^\infty \epsilon^{n + \frac{1}{m}} \phi_{mn + 1} \\
  + \ldots + 
  \sum_{n=n^*}^\infty \epsilon^{n + \frac{k}{m}} \phi_{mn + k} + \ldots
 +  \sum_{n=n^*}^\infty \epsilon^{n + \frac{m-1}{m}} \phi_{mn + m -1}.
\end{multline}
The bracketed terms are the early orders. Our interest is in deriving the form of the high-order terms, given by $\phi_{nm+k}$ for $k = 0, 1, \ldots, m-1$ as $n\to \infty$, which we claim can be represented using $m$ distinct ansatzes. 

However, the Saffman-Taylor problem turns out to be simpler because the late-order terms are only coupled a distance $\mathcal{O}(\epsilon)$ apart. In other words, for that problem, as $n \to \infty$, the leading-order behaviour of $\phi_{mn}$ in (\ref{eq:nhull_phiasym_split}) only depends on terms within its own sub-expansion. Because of this unique property, the late terms of the Saffman-Taylor problem are still given by a factorial over power ansatz. This simplification does not hold for our problem.

There exist other distinguished limits that may yield interesting results or alternative methodologies. For instance, what occurs if $\sigma_1$ and $\sigma_2$ in \eqref{qs_2hull} were to vary as $\epsilon \to 0$? Since the values of $\sigma_1$ and $\sigma_2$ determine the behaviour of the singulants, $\chi_k$ in \eqref{chik}, changing these powers changes the Stokes line configurations in the complex plane. There are indeed limitations of our methodology (also noted in \cite{trinh_2014}) where if the $\sigma_1$ and $\sigma_2$ are not rational numbers, then the asymptotics procedure becomes more difficult. Indeed, in some situations, it may be important to understand how $\sigma_1$ and $\sigma_2$ continuously approach some special value as $\epsilon \to 0$. For instance, Lustri \etal \cite{lustri_2012} study free-surface flow over an inclined step where the bottom topography is designed in order to produce leading-order downstream wave cancellation. However, as the geometry approached this optimzed configuration (by varying the inclination angle in the step), they noted an interesting behaviour where the Stokes line contribution from the two singularities (the stagnation point and corner in the step) rapidly oscillated in and out of phase. This variant of the close-singularities problem will present an interesting problem for future work. 

\ack We gratefully thank the referees for their insightful comments and many helpful suggestions. PHT thanks Lincoln College (Oxford) for generous support.

\appendix

\section{Relationship to the nonlinear equations of potential flow} \label{app:realfluid}

\begin{figure}[t] \centering
\includegraphics{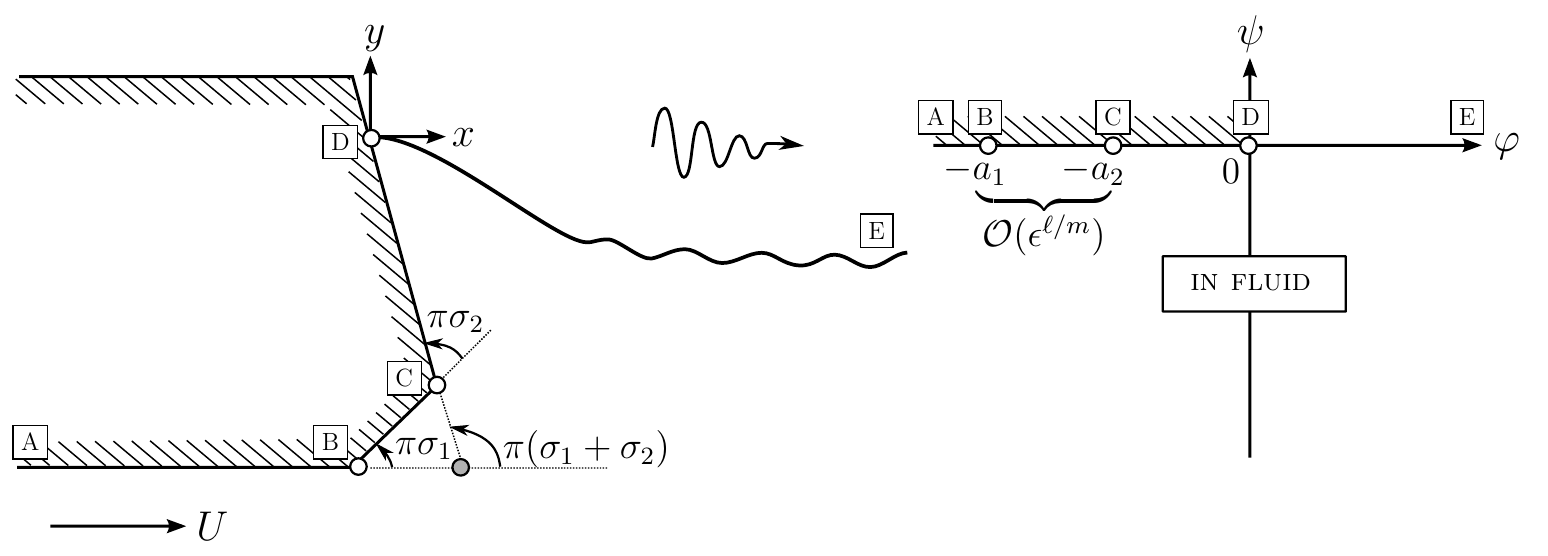}
\caption{Flow past a two-cornered ship in the (left) physical $xy$-plane and (right) complex potential $\varphi\psi$-plane. The two corners of divergent
angles $\pi\sigma_1$ and $\pi\sigma_2$ are mapped to $\varphi = -a_1, -a_2$ in the potential plane. The coalescing corners limit occurs if the two corners are $\Oh(\ep^{\ell/m})$ apart [see \eqref{eq:lmsig}], when the ship appears to possess only a single corner of angle $\pi(\sigma_1+\sigma_2)$. \label{fig:nhull}}
\end{figure}

Consider steady, incompressible, irrotational, and inviscid flow in the presence of gravity and past a semi-infinite body, which consists of a flat bottom and a piecewise linear front face, such as the one sketched in Figure \ref{fig:nhull}. The dimensional problem can be reposed in terms of a non-dimensional boundary-integral formulation in the potential $(\varphi, \psi)$-plane. The unknowns are the fluid speed $q = q(\varphi, \psi)$, and streamline angles, $\theta = \theta(\varphi, \psi)$, measured from the positive $x$-axis. The body and free-surface is given by the streamline $\psi = 0$, with $\psi \leq 0$ within the fluid. We assume the free-surface ($\varphi > 0$) attaches to the hull ($\varphi < 0$) at a stagnation point. The free-surface, with $\psi = 0$, is then obtained by solving a boundary-integral equation, coupled with Bernoulli's condition:
\begin{subequations} \label{governingreal}
\begin{gather}
 \log{q} = \log q_s + \frac{1}{\pi} \dashint_{0}^{\infty}
\frac{\theta(t)}{t - \varphi} \ \de{t} \label{eq:nhull_bdint} \\
 \epsilon q^2 \dd{q}{\varphi} = -\sin{\theta}, \quad \text{on $\varphi = 0$,}
\label{eq:nhull_bern}
\end{gather}
\end{subequations}
where $\epsilon = U^2/(gL)$ is the square of the Froude draft number with upstream flow $U$, gravity $g$, length scale $L = K/U$, and potential scale $K$. The function $q_s = q_s(\varphi)$ is calculated through 
\begin{equation}
 \log q_s \equiv \frac{1}{\pi} \dashint_{-\infty}^{0}
\frac{\theta(t)}{t - \varphi} \ \de{t},
\end{equation}
and for a surface-piercing object such as a ship, it is assumed that this function is known through the specification of the angle, $\theta$, along the negative real $\varphi$-axis, which corresponds to the geometry of the hull. For instance, a semi-infinite ship hull with a single corner at $\varphi = -a$ and with a face of angle $\pi \sigma$ is given by the $q_s$ in \eqref{qs_one}. Similarly, the two-singularity $q_s$ in \eqref{qs_2hull} corresponds to a two-cornered hull with divergent angles $\pi \sigma_1$ and $\pi\sigma_2$. See \cite{trinh_2011, trinh_2014} for more details on the derivation of (\ref{eq:nhull_bdint}) and (\ref{eq:nhull_bern}) as it pertains to the ship problem. For a more detailed reference on boundary integral equations for free-surface flows, see Vanden-Broeck \cite{vanden-broeck_2010}. 


In the exponential asymptotics framework, we are interested in studying the analytic continuation of the free-surface and thus allowing $q(\varphi, 0) \mapsto q(w)$ and $\theta(\varphi, 0) \mapsto \theta(w)$ to be complex functions of the complex potential, $\varphi + i0 \mapsto w$. However, in previous exponential asymptotic studies of the boundary integral equations of potential flow, it has been shown that because the exponential switched-on across Stokes line is almost exclusively determined by the local behaviour of the solution near the singularities in the complex plane (far away from where the boundary integral is evaluated), then the integral can be dropped from the analysis with little consequence to the methodology. This simplification is discussed in more detail within \cite{trinh_2011}, and is argued rigorously for the related problem of Saffman-Taylor viscous fingering in \cite{xie_2002}. 

In fact, by following the full procedure in \cite{trinh_2011}, it can be verified that the integral term serves to only change the amplitude coefficient of the waves [\eg in formula \eqref{phiexp_abs}] by a non-zero $\Oh(1)$ factor. This way, we simplify the full problem in \eqref{governingreal} to a simpler nonlinear differential equation in $q$. To be specific, analytic continuation of the integral equation \eqref{eq:nhull_bdint} into the upper-half plane yields
\begin{equation}
\log q + \im \theta = \log q_s + \frac{1}{\pi} \int_0^\infty \frac{\theta(t)}{t - w} \, \de{t},
\end{equation}
where note the $\im\theta$ term disappears and we recover the principal value integral \eqref{eq:nhull_bdint} by taking $w$ to the $\varphi$ axis from the upper-half plane. Now using the replacement $\im\theta \sim \log (q_s/q)$ instead of \eqref{eq:nhull_bdint}, and removing $\theta$ from Bernoulli's equation, yields the single differential equation
\begin{equation} \label{eq:nhull_simpbern}
\epsilon q_s q^3 \dd{q}{w} + \frac{\im}{2} \biggl[q^2 - q_s^2\biggr] = 0, \qquad q(0) = 0,
\end{equation}
where now $q$ is a complex-valued function, and the physical fluid boundary is identified with $w \in \mathbb{R}^+$. The boundary condition $q(0) = 0$ imposes a stagnation point attachment between the free-surface and solid boundary. Substitution of $\phi = q^2$ (note $\phi$ is not related to the fluid potential, $\varphi$, introduced earlier) leads to the main differential equation of this paper \eqref{phi_ode}. Note the qualitative similarities between real and complex parts of solutions to \eqref{eq:nhull_simpbern} shown in Figure~\ref{phiprofile}, and the full boundary-integral solutions displayed in Figure~7 of \cite{trinh_2011} and Figure~4 of \cite{trinh_2014}.

\section{General methodology for the close-singularity analysis} \label{sec:general}

In the paper, we principally studied the differential equation \eqref{phi_ode} subject to the two-singularity forcing, $q_s$, in \eqref{qs_2hull} where $\sigma_1 + \sigma_2 = 1/3$. These values were chosen so that the asymptotic expansion of $q_s$ (in powers of $\ep^{1/2}$) interleaves in the simplest non-trivial way with the usual $\epsilon$ expansion from the differential equation \eqref{phi_ode}. Other combinations of $\sigma_1 + \sigma_2$ will produce more complicated interleaving behaviour, and hence more complicated exponential-over-power divergence. In this section, we present the main ideas as to how the methodology is extended to the general problem. 

Rather than beginning with the outer analysis (far from singularities) as we have in Sec.~\ref{sec:cc13}, we will begin instead with the inner analysis, since it is easier to see the emergence of the correct divergent ansatz in this regime. 

\subsection{Step 1: Inner analysis with $w + a = \mathcal{O}(\epsilon^{\ell/m})$}

Our plan is to derive the leading-order inner problem, and then obtain the correct form of the solution as we tend to the outer region. Writing $q_0 = [w/(w+a)]^\sigma$ for the merged forcing function in \eqref{eq:qspre}, we note that as $w \to -a$, $q_s \sim c (w + a)^{-\sigma}$, where $c = (-a)^\sigma$, and $\sigma = \sigma_1 + \sigma_2$. The inner variables, $z$ and $\widehat{\phi}$, are defined using
\begin{equation}
w + a = (\epsilon/X)^{\ell/m} z^\ell, \qquad 
\phi = c^2(w+a)^{-2\sigma} \widehat{\phi}.	
\end{equation}
where we have written $X = \im/[c^3 (1 + 3\sigma)]$. The constants $c$ and $X$ are introduced to simplify the algebra. From (\ref{phi_ode}), the inner equation then becomes
\begin{equation} \label{eq:nhull_cc_innereq}
\frac{1}{m z^m} \Bigl[ -2\sigma \ell \widehat{q_s} \widehat{\phi}^2 + z \widehat{q_s} \widehat{\phi} \dd{\widehat{\phi}}{z} \Bigr] + \Bigl[\widehat{\phi} - \widehat{q_s}^2 \Bigr]
= 0. 
\end{equation}
In \eqref{eq:nhull_cc_innereq}, we have written $q_s = q_0 \widehat{q_s}$, and from \eqref{eq:nhull_cc_qs},
\begin{subequations}
\begin{equation} \label{qshat_appendix}
\widehat{q}_s = \frac{(X^{-\ell/m}z^\ell)^{\sigma_1+\sigma_2}}{(z^\ell + \beta)^{\sigma_1} (z^\ell - \beta)^{\sigma_2}} \equiv \sum_{n=0}^\infty \frac{\widehat{e_n}}{z^n},	
\end{equation}
where we have defined
\begin{equation}
\widehat{e_n} = 
\begin{cases}
\beta^{n/\ell} X^{n/m} f\left(\frac{n}{\ell}\right) & 
\text{\quad if $\text{mod}(n,\ell) = 0$} \\ 
0 & \text{\quad if $\text{mod}(n,\ell) \neq 0$.} 
\end{cases}	
\end{equation}

\end{subequations}

We wish to study the leading-order solution, $\widehat{\phi}$, as it tends outwards. Thus as $z \to \infty$, we substitute $\widehat{\phi} = \sum_{n=0}^\infty A_n/z^n$ into the inner equation (\ref{eq:nhull_cc_innereq}), giving
\begin{subequations}
\begin{gather}
A_0 = 1, \qquad A_n = \sum_{j=0}^n \widehat{e}_j\widehat{e}_{n-j} \quad \text{for $n < m$,} \\
 A_n = \sum_{j=0}^n \widehat{e}_j \widehat{e}_{n-j} + 
 \sum_{k=0}^{n-m} \widehat{e}_k \Biggl[ \sum_{j=0}^{n-m-k} \left( \frac{j + 2\sigma \ell}{m} \right)  A_j A_{n-m-k-j} \Biggr]. \qquad \text{for $n \geq m$}, \label{eq:nhull_cc_fullrr}
\end{gather}
\end{subequations}

A simple numerical computation assures us that $A_n$ diverges in the limit $n\to \infty$, and indeed, the form of the late terms follows
\begin{equation} \label{nhull_cc_Aninner}
 A_n \sim \Omcc \e^{\im \tau} \Gamma\left( \frac{n}{m} + \gamma \right)
\exp \left[ \sum_{j=1}^{m-1} \mu_j n^{\frac{m-j}{m}}\right],
\end{equation}

\noindent as $n \to \infty$, with $\gamma, \mu_j \in \mathbb{C}$, $\Omcc, \tau \in \mathbb {R}$. Like the discussion surrounding \eqref{An_in_both}, we caution the reader that the form of \eqref{nhull_cc_Aninner} may need to be multipled by $(-1)^n$ depending on the form of \eqref{nhull_cc_Aninner}. The special case of $\mu_j \equiv 0$ for all $j$ corresponds to the more typical case of factorial divergence, which is observed for most problems in exponential asymptotics. The suitability of (\ref{nhull_cc_Aninner}) can be understood by dividing (\ref{eq:nhull_cc_fullrr}) by $A_n$ and writing it as
\begin{equation} \label{nhull_cc_Aninnerb}
 1 + 
 \sum_{j=0}^m 
 \Biggl[\sum_{k=0}^j \Bigl(-\frac{1}{m}\Bigr)\widehat{e_k} A_{j-k} \Biggr]
  \frac{n A_{n-m-j}}{A_n} + 
  \Biggl[\frac{(4\ell-m)\widehat{e_0} A_0}{3m}\Biggr] \frac{A_{n-m}}{A_n} +
\ldots \ = 0
\end{equation}

\noindent Then using the ansatz (\ref{nhull_cc_Aninner}), we can see that the series expansions of the various ratios appearing in (\ref{nhull_cc_Aninnerb}) are given by
\begin{center}
\begin{tabular}{ccccccccc}
$\dfrac{A_{n-m}}{A_n}$ 
&$\sim$
& $\dfrac{c_{10}}{n}$ 
&$+$ 
& $\dfrac{c_{11}}{n^{\frac{m+1}{m}}}$   
& $+\quad \ldots \quad +$
& $\dfrac{c_{1(m-1)}}{n^{\frac{2m-1}{m}}}$ 
&$+$ 
& $\dfrac{c_{1m}}{n^2}$ \\
$\dfrac{A_{n-m-1}}{A_n}$ 
&$\sim$
& 
& 
& $\dfrac{c_{21}}{n^{\frac{m+1}{m}}}$   
&$+ \quad \ldots \quad +$
& $\dfrac{c_{2(m-1)}}{n^{\frac{2m-1}{m}}}$ 
&$+$ 
& $\dfrac{c_{2m}}{n^2}$ \\
$\vdots$ 
&
&
&
&
& $\ddots$
&
& 
& $\vdots$ \\
$\dfrac{A_{n-2m+1}}{A_n}$ 
&$\sim$
&
&
& 
& $\ldots$ 
& $\dfrac{c_{(m-1)(m-1)}}{n^{\frac{2m-1}{m}}}$ 
& $+$ 
& $\dfrac{c_{(m-1)m}}{n^2}$ \\
$\dfrac{A_{n-2m}}{A_n}$ 
&$\sim$
&
&
&
& $\ldots$ 
&
& 
&$\dfrac{c_{mm}}{n^2}$, \\
\end{tabular}
\end{center}

\noindent where the factors $c_{jk}$ are functions of $\mu_j$ and $\gamma$ and
in general, require the higher-order terms of Stirling's approximation to the
Gamma function, and can be derived using a computer algebra system.

Thus, when we search for the $n\to \infty$ limit in (\ref{eq:nhull_cc_fullrr}), we will
need to conserve the terms with factors of $A_n$ and $A_{n-m}$, and
terms of orders $nA_{n-m}$ to $nA_{n-2m}$. Terms with other indices will be of lower
order. Then, for equation (\ref{nhull_cc_Aninnerb}) and using the ansatz, the leading order behaviour at $\mathcal{O}(n^{-1})$ is automatically satisfied by the Gamma function; $\mu_1$ is then determined at $\mathcal{O}(n^{-\frac{m+1}{m}})$, $\mu_2$ is determined at $\mathcal{O}(n^{-\frac{m+2}{m}})$, and so on until $\mu_{m-1}$ is determined at $\mathcal{O}(n^{-\frac{2m-1)}{m}})$. Lastly, the constant $\gamma$ is determined at $\mathcal{O}(n^{-2})$. The prefactor $\Omcc \e^{\im\tau}$ can then be
found by numerically solving the recurrence relation in (\ref{eq:nhull_cc_fullrr}). 

\subsection{Step 2: Outer analysis with $w + a = \mathcal{O}(1)$}

We are now in a position to return to the outer equation of (\ref{phi_ode}). Upon substituting the expansion $\phi = \sum \epsilon^{n/m} \phi_n$ into the equation, the first $m$ terms are
\begin{equation}
\phi_0 = q_0^2, \qquad \phi_n = q_0^2\sum_{k=0}^n e_{k} e_{n-k} \qquad \text{for $n <
m$,}	
\end{equation}
and motivated by the form of the inner ansatz (\ref{nhull_cc_Aninner}), we assume that as $n \to \infty$, the late terms are given by
\begin{equation}
\phi_{n} \sim \frac{P(w) \Gamma\left( \frac{n}{m} + \gamma\right)}{[\chi(w)]^{\frac{n}{m}+\gamma}}
\exp\left[ \sum_{i=1}^{m-1} r_i(w)
n^{\frac{m-i}{m}}\right],	
\end{equation}
for complex functions $\chi$, $r_i$, and $P$, and constant $\gamma$. Substitution of this ansatz into the differential equation yields a similar procedure to that which we used for the inner equation, except that now, we wish to keep terms $\phi_n$ and $\phi_{n-m}$, as well as derivatives $\phi'_{n-m}$ to $\phi'_{n-2m}$. The relevant terms at $\mathcal{O}(\epsilon^{\frac{n}{m}})$ are then
\begin{equation} \label{eq:nhull_cc_genorderphi}
 \phi_n + 
 \sum_{j=0}^m 
 \Biggl[\sum_{k=0}^j \Bigl(-\im q_0\Bigr)e_k \phi_{j-k} \Biggr] \phi'_{n-m-j} + 
  \Biggl[(-\im) e_0 \phi'_0 \Biggr]\phi_{n-m} + \ldots \ = 0.
\end{equation}

As $n\to\infty$, the leading-order terms involve $\phi_n + (-\im q_0) e_0 \phi_0 \phi'_{n-m} = 0$, with $e_0 = 1$ and $\phi_0 = q_0^2$. Thus the equation for $\chi$ is 
\begin{equation}
\chi(w) = \int_{-a}^w \frac{\im}{q_0^3(\varphi)} \ d\varphi,	
\end{equation}
The subsequent orders each yield first-order differential equations for $r_i(w)$ and $P(w)$. In order to match with the form of the inner ansatz in \eqref{nhull_cc_Aninner}, we require $r_i(-a) = \mu_i$. Once the boundary conditions on $r_i(w)$ are imposed, the value of $\gamma$ can be verified by setting $n = 0$ in the late-orders ansatz and matching its behaviour as $w \to -a$. This will turn out to be the same $\gamma$ computed directly from the inner analysis. 

\subsection{Step 3: Optimal truncation and Stokes smoothing}

Once the late-order terms, $\phi_n$, have been found, we can re-scale near the Stokes lines and derive the exponential switchings. We begin by truncating the asymptotic expansion at $n = \mathcal{N}$,
\begin{equation}
\phi = \sum_{n=0}^{\mathcal{N}-1} \epsilon^{\frac{n}{m}} \phi_n + R_\mathcal{N}	
\end{equation}
and substitute this expression into \eqref{phi_ode}. This gives a linear equation for the remainder, $R_\mathcal{N}$, as well as a single $\mathcal{O}(\epsilon^{\mathcal{N}/m})$ term:
\begin{equation} \label{eq:nhull_cc_remeq}
 \L(R_\mathcal{N}; \epsilon) \sim \epsilon^{\mathcal{N}/m}
\phi_\mathcal{N},
\end{equation}

\noindent where the relevant terms of the linear operator $\L$ are given by
\begin{equation} \label{eq:nhull_cc_Lop}
 \L(R_\mathcal{N}; \epsilon) = R_\mathcal{N} + 
 \sum_{j=0}^m 
 \Biggl[\sum_{k=0}^j \Bigl(-\im q_0\Bigr)e_k \phi_{j-k} \Biggr]
\epsilon^{1+\frac{j}{m}} \dd{R_\mathcal{N}}{w} + 
  \Biggl[(-\im) e_0 \phi'_0 \Biggr]\epsilon R_\mathcal{N} + \ldots \ = 0.	
\end{equation}

\noindent Examine now Table \ref{tab:orders}, which summarises the connection between the inner, outer, and Stokes Smoothing procedures. To review: the table shows that for the inner analysis, seeking the $n\to \infty$ limit for $A_n$ involves equations in orders of $n^{-1/m}$; for the outer analysis, the $n\to \infty$ limit also involves equations in orders of $n^{-1/m}$; for the Stokes smoothing procedure, the exponential behaviour of $R_N$ is determined by matching at orders in $\epsilon^{1/m}$. The equations at each order are not the same for all three analyses, but they share corresponding values near the inner region (between terms $A_n$ and $\phi_n$) and near the Stokes line (between terms $\phi_n$ and $R_N$).

\newcolumntype{Z}{>{\centering\arraybackslash}X}
\newcolumntype{Y}{>{\raggedright\arraybackslash}X}
\newcolumntype{W}{>{\raggedleft\arraybackslash}X}
\renewcommand{\tabularxcolumn}[1]{>{\arraybackslash}m{#1}}
\newcommand{\otoprule}{\midrule[0.3pt]}
\begin{table} \centering
\begin{tabularx}{0.98\textwidth}{YYZZZZZZ} \toprule
\textsc{region} & \textsc{form} & \textsc{first} & \textsc{first} &
\textsc{second} & \ldots & \textsc{last} & \textsc{last} \\ \otoprule 
Inner  & $\displaystyle \sum_n \frac{A_n}{z^n}$ & $A_n$ & $nA_{n-m}$ &
$nA_{n-m-1}$ & \ldots & $nA_{n-2m}$ & $A_{n-m}$ \\ \\[0.1cm]
Outer  & $\displaystyle \sum_n \epsilon^{\frac{n}{m}} \phi_n$ & $\phi_n$ &
$\phi'_{n-m}$ & $\phi'_{n-m-1}$ & \ldots & $\phi'_{n-2m}$ & $\phi_{n-m}$ \\
\\[0.1cm]
Stokes & $\displaystyle R_\mathcal{N} e^{F}$ & $R_\mathcal{N}$ &
$\epsilon^{\frac{1}{m}} R'_\mathcal{N}$ & $\epsilon^{\frac{2}{m}}
R'_\mathcal{N}$ & \ldots & $\epsilon^{\frac{m-1}{m}} R'_\mathcal{N}$ & $\epsilon
R'_\mathcal{N}$ \\ \bottomrule 
\end{tabularx}
\caption{Connection between the expansion forms of the inner, outer, and Stokes line analyses. The columns indicate which terms from \eqref{nhull_cc_Aninnerb}, \eqref{eq:nhull_cc_genorderphi}, and \eqref{eq:nhull_cc_Lop} contribute and at which order during each of the respective analyses. For instance, the $\Gamma(n/m+\gamma)$ behaviour of the inner analysis is determined by a balance of terms $A_n$ and $nA_{n-m}$. \label{tab:orders}}
\end{table}

Before we solve the inhomogeneous equation (\ref{eq:nhull_cc_remeq}), we first seek the homogeneous solutions of $\L = 0$. We expand the remainder as
\begin{equation} \label{eq:nhull_cc_Rnhomo}
 R_\mathcal{N} = \overline{P}(w)\exp\Bigl[F(w)\Bigr] = 
 \overline{P}(w)\exp\left[\sum_{j=0}^{m-1}
\frac{F_i(w)}{\epsilon^{\frac{j}{m}}} \right].
\end{equation}
%
%
Substitution into $\L = 0$ yields equations at $\mathcal{O}(1)$, $\mathcal{O}(\epsilon^{\frac{1}{m}})$, $\ldots$, $\mathcal{O}(\epsilon^{\frac{m-1}{m}})$, $\mathcal{O}(\epsilon)$, which determines $F'_0$, $F'_1$, $\ldots$ $F'_{m-1}$, and $\overline{P}$, each in turn.

To solve the inhomogeneous equation (\ref{eq:nhull_cc_remeq}), we multiply the homogeneous solution (\ref{eq:nhull_cc_Rnhomo}) by the Stokes smoothing factor $\mathcal{S}$, then substitute $R_\mathcal{N} = \mathcal{S}\overline{P}\exp[F(w)]$ this expression into (\ref{eq:nhull_cc_remeq}), giving
\begin{equation} \label{finalRNeq}
-\epsilon \im q_0^3 \dd{\mathcal{S}}{w} \overline{P}\e^{F} \sim \epsilon^{\mathcal{N}/m} 
\phi_\mathcal{N}, 	
\end{equation}
where the bracketed factor corresponds to the $j = 0$ contribution from the summation in (\ref{eq:nhull_cc_Lop}). Finally, we can establish a relationship between the components of the final exponential given by $R_\mathcal{N}$, and the numerical constant of the inner problem $\Omcc$, all related through the late-order terms $r_i$ and $P$. Simplification of \eqref{finalRNeq} will then provide an equation for the value of $\mathcal{S}$ across the Stokes line.

\bibliographystyle{plainnat}
\bibliography{mybib}

\end{document}